\documentclass{amsart}
\usepackage[utf8]{inputenc}
\pdfoutput=1

\usepackage{amsmath,amsfonts,amssymb}
\usepackage{amsthm}
\usepackage{mathtools}
\usepackage{multirow}
\usepackage{tikz}
\usetikzlibrary{calc}
\usetikzlibrary{shapes.geometric}
\usetikzlibrary{arrows}
\usetikzlibrary{automata}
\usetikzlibrary{positioning}
\usepackage{url}
\usepackage{hyperref}
\usepackage{doi}
\usepackage{xcolor}
\usepackage{booktabs}
\usepackage{xspace}
\usepackage[all,cmtip]{xy}
\usepackage{xargs}

\theoremstyle{plain}

\newtheorem{definition}{Definition}
\theoremstyle{remark}
\newtheorem{remark}{Remark}

\newcommand{\PP}{\ensuremath{\mathbb{P}}}
\newcommand{\FF}{\ensuremath{\mathbb{F}}}
\newcommand{\FFbar}{\ensuremath{\overline{\mathbb{F}}}}
\newcommand{\ZZ}{\ensuremath{\mathbb{Z}}}
\newcommand{\XC}{\ensuremath{\mathcal{C}}}
\newcommand{\EC}{\ensuremath{\mathcal{E}}}
\newcommand{\AV}{\ensuremath{\mathcal{A}}}

\newcommand{\field}{\ensuremath{\Bbbk}\xspace}
\newcommand{\fieldbar}{\ensuremath{\overline{\Bbbk}}\xspace}
\newcommand{\Jac}[1]{\ensuremath{\mathcal{J}({#1})}}
\newcommand{\Aut}{\ensuremath{\mathrm{Aut}}}
\newcommand{\RAut}{\ensuremath{\mathrm{RA}}}

\newcommand{\dualof}[1]{\ensuremath{#1^\dagger}}

\newcommand{\classof}[1]{\ensuremath{\big[{#1}\big]}}

\newcommand{\subgrp}[1]{\ensuremath{\langle{#1}\rangle}}

\newcommandx{\IG}[3][1=g, 2=\ell, 3=p]{\ensuremath{\Gamma_{#1}(#2;#3)}\xspace}
\newcommandx{\SSIG}[3][1=g, 2=\ell, 3=p]{\ensuremath{\Gamma^{SS}_{#1}(#2;#3)}\xspace}
\newcommand{\RichelotIG}{\ensuremath{\Gamma^{SS}_{2}(2;p)}\xspace}

\newcommand{\TypeA}{\textsf{Type-A}\xspace}
\newcommand{\TypeI}{\textsf{Type-I}\xspace}

\newcommand{\TypeII}{\textsf{Type-II}\xspace}
\newcommand{\TypeIII}{\textsf{Type-III}\xspace}

\newcommand{\TypeIV}{\textsf{Type-IV}\xspace}

\newcommand{\TypeV}{\textsf{Type-V}\xspace}
\newcommand{\TypeVI}{\textsf{Type-VI}\xspace}

\newcommand{\TypeExE}{{\textsf{Type-$\Pi$}}\xspace}
\newcommand{\TypeEsquared}{{\textsf{Type-$\Sigma$}}\xspace}
\newcommand{\TypeExEzero}{{\textsf{Type-$\Pi_0$}}\xspace}
\newcommand{\TypeEzerosquared}{{\textsf{Type-$\Sigma_0$}}\xspace}
\newcommand{\TypeExEtwelvecubed}{{\textsf{Type-$\Pi_{12^3}$}}\xspace}
\newcommand{\TypeEtwelvecubedsquared}{{\textsf{Type-$\Sigma_{12^3}$}}\xspace}
\newcommand{\TypeEzeroxEtwelvecubed}{{\textsf{Type-$\Pi_{0,12^3}$}}\xspace}

\newcommand{\TypeExElabel}{\ensuremath{\Pi}}
\newcommand{\TypeExEzerolabel}{\ensuremath{\Pi_0}}
\newcommand{\TypeExEtwelvecubedlabel}{\ensuremath{\Pi_{12^3}}}
\newcommand{\TypeEzeroxEtwelvecubedlabel}{\ensuremath{\Pi_{0,12^3}}}
\newcommand{\TypeEsquaredlabel}{\ensuremath{\Sigma}}
\newcommand{\TypeEzerosquaredlabel}{\ensuremath{\Sigma_0}}
\newcommand{\TypeEtwelvecubedsquaredlabel}{\ensuremath{\Sigma_{12^3}}}
\newcommand{\TypePhilabel}{\ensuremath{\Phi}}

\title{An atlas of the Richelot isogeny graph}

\author{Enric Florit}
\address{IMUB - Universitat de Barcelona, Gran Via de les Corts Catalanes 585, 08007 Barcelona, Spain}
\email{efz1005@gmail.com}

\author{Benjamin Smith}
\address{Inria and Laboratoire d'Informatique de l'École polytechnique (LIX), Institut Polytechnique de Paris, 1 rue Honoré d'Estienne d'Orves, 91120 Palaiseau, France}
\email{smith@lix.polytechnique.fr}

\thanks{The second author was supported by ANR CIAO}

\date{December 2020}

\begin{document}

\begin{abstract}
    We describe and illustrate
    the local neighbourhoods of vertices and edges
    in the \((2,2)\)-isogeny graph
    of principally polarized abelian surfaces,
    considering the action of automorphisms.
    Our diagrams are intended to build intuition
    for number theorists and cryptographers
    investigating isogeny graphs in dimension/genus 2,
    and the superspecial isogeny graph in particular.
\end{abstract}

\maketitle

\section{
    Introduction
}

This article is an illustrated guide to the Richelot isogeny graph.
Following Katsura and Takashima~\cite{2020/Katsura--Takashima},
we present diagrams of the neighbourhoods of general vertices of each type.
Going further,
we also compute diagrams of neighbourhoods of general edges,
which can be used to glue the vertex neighbourhoods together.
Our aim is to build intuition on the various combinatorial
structures in the graph,
providing concrete examples for some of the more pathological cases.
The authors have used the results presented here
to verify computations and form conjectures
when investigating the behaviour of random walks
in superspecial isogeny graphs~\cite{2020/Florit--Smith}.

We work over a ground field \field
of characteristic not \(2\), \(3\), or \(5\).
In our application to superspecial PPASes,
\(\field = \FF_{p^2}\),
though our computations were mostly done over function fields over cyclotomic
fields.

Let \(\AV/\field\) be a principally polarized abelian surface (PPAS).
A \emph{\((2,2)\)-isogeny}, or \emph{Richelot isogeny},
is an isogeny \(\phi: \AV \to \AV'\) of PPASes
whose kernel is a maximal \(2\)-Weil isotropic subgroup of \(\AV[2]\).
Such a \(\phi\)
has kernel isomorphic to \((\ZZ/2\ZZ)^2\);
it respects the principal polarizations
\(\lambda\) and \(\lambda'\) on \(\AV\) and \(\AV'\), respectively,
in the sense that
\(\phi^*(\lambda') = 2\lambda\);
and its (Rosati) \emph{dual isogeny} \(\dualof{\phi}: \AV' \to \AV\)
satisfies \(\dualof{\phi}\circ\phi = [2]_{\AV}\).

The \emph{\((2,2)\)-isogeny} or \emph{Richelot isogeny graph} 
is the directed weighted multigraph defined as follows.
The \emph{vertices} are isomorphism classes of
PPASes over \(\field\).
If \(\AV\) is a PPAS,
then \(\classof{\AV}\) denotes the corresponding vertex.
The \emph{edges} are isomorphism classes of
\((2,2)\)-isogenies
(\(\phi_1: \AV_1 \to \AV_1'\) 
and \(\phi_2: \AV_2 \to \AV_2'\) 
are isomorphic
if there are isomorphisms
of PPASes
$\alpha \colon \AV_1 \to \AV_2$
and
$\beta\colon \AV_1' \to \AV_2'$ 
such that 
\(\phi_2\circ\alpha = \beta\circ\phi_1\)).

The edges are \emph{weighted} by the number of distinct kernels yielding
isogenies in their class.
The weight of an edge \(\classof{\phi}\)
is denoted by \(w(\classof{\phi})\).
If \(\classof{\phi}: \classof{\AV}\to \classof{\AV'}\) is an edge,
then $w(\classof{\phi}) = n$
if and only if 
there are \(n\) kernel subgroups \(K \subset \AV[2]\) 
such that \(\AV' \cong \AV/K\)
(this is independent of the choice of
representative isogeny \(\phi\)).

There are fifteen maximal \(2\)-Weil-istropic subgroups in \(\AV[2]\),
though some (or all) might not be defined over \(\field\).
The sum of the weights of the edges leaving any vertex
is therefore at most 15.

The isogeny graph breaks up into connected components within isogeny classes.
We are particularly interested in the superspecial isogeny class.
Recall that a PPAS \(\AV/\FFbar_p\) is \emph{superspecial}
if its Hasse--Witt matrix vanishes identically.
Equivalently, \(\AV\) is superspecial
if it is isomorphic \emph{as an unpolarized abelian variety}
to a product of supersingular elliptic curves.
For background on superspecial and supersingular
abelian varieties in low dimension, 
we refer to 
Ibuyiyama, Katsura, and Oort~\cite{1986/Ibukiyama--Katsura--Oort}
and Brock's thesis~\cite{1993/Brock}.
For more general results, we refer to Li and Oort~\cite{1998/Li--Oort}.

\begin{definition}
    The superspecial Richelot isogeny graph
    is the subgraph \(\RichelotIG\) of 
    the Richelot isogeny graph over \(\FF_{p^2}\)
    supported on the superspecial vertices.
\end{definition}

Recall that \RichelotIG has
\(p^3/2880 + O(p)\) vertices
(see~\S\ref{sec:superspecial-count} for a more precise statement),
and is connected~\cite{2020/Jordan--Zaytman}.
If \(\AV/\FFbar_p\) represents a vertex in \RichelotIG,
then the invariants corresponding to \(\classof{\AV}\)
are defined over \(\FF_{p^2}\),
as are all 15 of the \((2,2)\)-isogeny kernels---so
\RichelotIG is a \(15\)-regular graph.
It has interesting number-theoretic properties and applications
(such as Mestre's \emph{méthode des graphes}~\cite{1986/Mestre}),
and potential cryptographic applications
(including~\cite{2009/CGLgenusg,
2018/Takashima,
2019/Flynn--Ti,
2020/Castryck--Decru--Smith,
2020/Costello--Smith}).
All of these applications
depend on a clear understanding of the structure of~\RichelotIG:
for example,
the local neighbourhoods of vertices with extra automorphisms
(and their inter-relations)
affect
the expansion properties
and random-walk behaviour of \RichelotIG,
as we see in~\cite{2020/Florit--Smith}.

\section{
    Richelot isogenies and isogeny graphs
}
\label{sec:basics}

There are two kinds of PPASes:
products of elliptic curves (with the product polarization)
and Jacobians of genus-2 curves.
The algorithmic construction of isogenies 
depends fundamentally on whether the PPASes are Jacobians or elliptic products.
We recall the Jacobian case 
in~\S\ref{sec:Richelot},
and the elliptic product case
in~\S\ref{sec:elliptic-product-isogenies}.

\subsection{Richelot isogenies}
\label{sec:Richelot}

Let \(\XC: y^2 = F(x)\) be a genus-2 curve,
with \(F\) squarefree of degree \(5\) or \(6\).
The kernels of \((2,2)\)-isogenies from \(\Jac{\XC}\)
correspond to factorizations of \(F\) into quadratics
(of which one may be linear, if \(\deg(F) = 5\)):
\[
    \XC: y^2 = F(x) = F_1(x)F_2(x)F_3(x)
    \,,
\]
up to permutation of the \(F_i\) and constant multiples.
We call such factorizations \emph{quadratic splittings}.
The kernel (and isogeny) is defined over \(\field\)
if the splitting is.

Fix one such 
quadratic splitting \(\{F_1,F_2,F_3\}\);
then the corresponding subgroup \(K\subset\Jac{\XC}[2]\) 
is the kernel of a \((2,2)\)-isogeny
\(\phi: \Jac{\XC} \to \Jac{\XC}/K\).
For each \(1 \le i \le 3\),
we write \(F_i(x) = F_{i,2}x^2 + F_{i,1}x + F_{i,0}\).
Now let
\[
    \delta 
    = 
    \delta(F_1,F_2,F_3)
    := 
    \begin{vmatrix}
        F_{1,0} & F_{1,1} & F_{1,2} 
        \\
        F_{2,0} & F_{2,1} & F_{2,2} 
        \\
        F_{3,0} & F_{3,1} & F_{3,2} 
    \end{vmatrix}
    \,.
\]

If \(\delta(F_1,F_2,F_3) \not= 0\),
then 
\(\Jac{\XC}/K\) is isomorphic to a Jacobian \(\Jac{\XC'}\),
which we can compute using Richelot's algorithm
(see~\cite{1988/Bost--Mestre} and~\cite[\S8]{2005/Smith}):
\(\XC'\) is defined by
\[
    \XC': y^2 = G_1(x)G_2(x)G_3(x)
    \quad
    \text{where}
    \quad
    G_i(x) := \frac{1}{\delta}(F_j'(x)F_k(x) - F_k'(x)F_j(x))
\]
for each cyclic permutation \((i,j,k)\) of \((1,2,3)\).
The quadratic splitting \(\{G_1,G_2,G_3\}\)
corresponds to the kernel of the dual isogeny 
\(\dualof{\phi}: \Jac{\XC'}\to\Jac{\XC}\).

If \(\delta(F_1,F_2,F_3) = 0\),
then \(\Jac{\XC}/K\) is isomorphic to 
an elliptic product \(\EC\times\EC'\).
which we can compute as follows.
There exist linear polynomials \(U\) and \(V\)
such that
\(F_1 = \alpha_1 U^2 + \beta_1 V^2\)
and \(F_2 = \alpha_2 U^2 + \beta_2 V^2\)
for some \(\alpha_1\), \(\beta_1\), \(\alpha_2\), and \(\beta_2\);
and since in this case 
\(F_3\) is a linear combination of \(F_1\) and \(F_2\),
we must have
\(F_3 = \alpha_3 U^2 + \beta_3 V^2\)
for some \(\alpha_3\) and \(\beta_3\).
The elliptic factors are defined by
\[
    \EC: y^2 = \prod_{i=1}^3(\alpha_i x + \beta_i)
    \qquad
    \text{and}
    \qquad
    \EC': y^2 = \prod_{i=1}^3(\beta_i x + \alpha_i)
    \,,
\]
and 
the isogeny \(\phi: \Jac{\XC} \to \EC\times\EC'\)
is induced by the product of the
double covers
\(\pi: \XC \to \EC\) resp. \(\pi': \XC \to \EC'\)
mapping \((x,y)\) to 
\((U^2/V^2,y/V^3)\) resp. \((V^2/U^2,y/U^3)\).

\subsection{Isogenies from elliptic products}
\label{sec:elliptic-product-isogenies}

Consider a generic pair of elliptic curves 
\begin{align*}
    \EC: y^2 & = (x - s_1)(x - s_2)(x - s_3)
    \,,
    \\
    \EC': y^2 & = (x - s_1')(x - s_2')(x - s_3')
    \,.
\end{align*}
We have
\( \EC[2] = \{ 0_{\EC}, P_1, P_2, P_3 \} \)
and
\( \EC'[2] = \{ 0_{\EC'}, P_1', P_2', P_3' \} \)
where \(P_i := (s_i,0)\)
and \(P_i' := (s_i',0)\).
For each \(1 \le i \le 3\),
we let 
\[
    \psi_i: \EC \longrightarrow \EC_i := \EC/\subgrp{P_i}
    \quad
    \text{and}
    \quad 
    \psi_i': \EC' \to \EC_i' := \EC'/\subgrp{P_i'}
\]
be the quotient \(2\)-isogenies.
These can be computed using Vélu's formul\ae{}~\cite{1971/Velu}.

Nine of the fifteen kernel subgroups of \((\EC\times\EC')[2]\)
correspond to products of elliptic \(2\)-isogeny kernels.
Namely, 
for each \( 1 \le i, j \le 3 \)
we have the kernel
\[
    K_{i,j}
    :=
    \subgrp{
        (P_i,0_{\EC'})
        ,
        (0_{\EC},P_i')
    }
    \subset
    (\EC\times\EC')[2]
\]
of the product isogeny
\[
    \phi_{i,j} 
    := 
    \psi_i\times\psi_j: 
    \EC\times\EC' 
    \longrightarrow 
    \EC_i\times\EC_j'
    \cong 
    (\EC\times\EC')/K_{i,j} 
    \,.
\]

The other six kernels correspond to \(2\)-Weil anti-isometries
\(\EC[2]\cong\EC'[2]\):
they are
\[
    K_{\pi} 
    := 
    \{ 
        (0_{\EC},0_{\EC'}),
        (P_{1},P_{\pi(1)}'),
        (P_{2},P_{\pi(2)}'),
        (P_{3},P_{\pi(3)}')
    \}
    \quad
    \text{for }
    \pi \in \operatorname{Sym}(\{1,2,3\})
    \,,
\]
with quotient isogenies
\[
    \phi_{\pi}: \EC\times\EC' \longrightarrow \AV_{\pi} := (\EC\times\EC')/K_{\pi}
    \,.
\]
If the anti-isometry \(P_i \mapsto P_{\pi(i)}'\)
is induced by an isomorphism \(\EC \cong \EC'\),
then \(\AV_{\pi}\cong\EC\times\EC'\).
Otherwise,
following~\cite[Prop.~4]{2000/Howe--Leprevost--Poonen},
\(\AV_{\pi}\) is the Jacobian of a genus-2 curve
\[
    \XC_{\pi}: y^2 = -F_1(x)F_2(x)F_3(x)
\]
where
\[
    F_i(x) := A(s_j-s_i)(s_i-s_k)x^2 
           + B(s'_{j}-s'_{i})(s'_{i}-s'_{k})
\]
for each cyclic permutation \((i,j,k)\) of \((1,2,3)\),
with
\begin{align*}
    A & := \frac{a_1}{a_2}\prod (s'_i - s'_j)^2 
    \text{ where }
    a_1 := \sum \frac{(s_j - s_i)^2}{s'_j - s'_i}
    \text{ and }
    a_2 := \sum s_i(s'_k - s'_j)
    \,,
    \\
    B & := \frac{b_1}{b_2}\prod (s_i - s_j)^2
    \text{ where }
    b_1 := \sum \frac{(s'_j - s'_i)^2}{s_j - s_i}
    \text{ and }
    b_2 := \sum s'_i(s_k - s_j)
    \,,
\end{align*}
where the sums and products
are over cyclic permutations \((i,j,k)\) of \((1,2,3)\).
The dual isogeny \(\dualof{\phi_{\pi}}: \Jac{\XC_{\pi}} \to
\EC\times\EC'\) corresponds 
to the splitting \(\{F_1,F_2,F_3\}\).

\section{
    Automorphism groups of abelian surfaces
}
\label{sec:ppas}

We now consider the impact of automorphisms on edge weights in the
isogeny graph,
following Katsura and Takashima~\cite{2020/Katsura--Takashima},
and recall the explicit classification
of reduced automorphism groups of PPASes.
In contrast with elliptic curves,
where (up to isomorphism) only two curves
have nontrivial reduced automorphism group,
with PPASes we see much richer structures
involving many more vertices.
Proofs for all of the results in this section
can be found
in~\cite{1986/Ibukiyama--Katsura--Oort},
\cite{2020/Katsura--Takashima},
and~\cite{2020/Florit--Smith}.

\subsection{Automorphisms and isogenies}

Let \(\phi : \AV \to \AV/K\) be a \((2,2)\)-isogeny
with kernel \(K\).
Let \(\alpha\) be an automorphism of \(\AV\),
and let \(\phi': \AV \to \AV/\alpha(K)\)
be the quotient isogeny;
then \(\alpha\)
induces an isomorphism
\(\alpha_*: \AV/K \to \AV/\alpha(K)\)
such that 
\(\alpha_*\circ\phi = \phi'\circ\alpha\).

If \(\alpha(K) = K\),
then \(\AV/K = \AV/\alpha(K)\),
so \(\alpha_*\) is an automorphism of \(\AV/K\).
Going further, if \(S\) is the stabiliser
of \(K\) in \(\Aut(\AV)\),
then \(S\) induces an isomorphic subgroup \(S'\) of \(\Aut(\AV/K)\),
and in fact \(S'\) is the stabiliser of \(\ker(\dualof{\phi})\) in
\(\Aut(\AV/K)\).

If \(\alpha(K) \not= K\) then
the quotients \(\AV/K\) and \(\AV/\alpha(K)\)
are different,
so \(\alpha_*\) is an isomorphism but not an automorphism.
The isogenies \(\phi\) and \(\phi_\alpha := \alpha_*^{-1}\circ\phi'\)
have identical domains and codomains, but distinct kernels;
thus, they both represent the same edge in the isogeny graph,
and \(w(\classof{\phi}) > 1\).

Every PPAS has \([-1]\) in its automorphism group,
but \([-1]\) fixes every kernel and commutes with
every isogeny---so it has no impact on edges or weights in the isogeny
graph.
We can therefore simplify by quotienting
\([-1]\) out of the picture.

\begin{definition}
    If \(\AV\) is a PPAS,
    then its \textbf{reduced automorphism group} is 
    \[
        \RAut(\AV) := \Aut(\AV)/\subgrp{[-1]}
        \,.
    \]
\end{definition}

Since \(\subgrp{[-1]}\) is contained in the centre of \(\Aut(\AV)\),
the quotient \(\RAut(\AV)\) acts on the set of kernel subgroups of
\(\AV[2]\).
We have two useful results
for \((2,2)\)-isogenies \(\phi: \AV \to \AV/K\).
First,
if \(O_K\) is the orbit of \(K\) under \(\RAut(\AV)\),
then there are \(\#O_K\) distinct kernels
of isogenies representing~\(\classof{\phi}\):
that is,
\[
    w(\classof{\phi}) = \#O_K\,.
\]
Second,
we have the ``ratio principle''
from~\cite[Lemma~1]{2020/Florit--Smith}:
\begin{equation}
    \label{eq:ratio-principle}
    \#\RAut(\AV)\cdot w(\classof{\dualof{\phi}})
    =
    \#\RAut(\AV')\cdot w(\classof{\phi})
    \,.
\end{equation}

\subsection{Reduced automorphism groups of Jacobians}

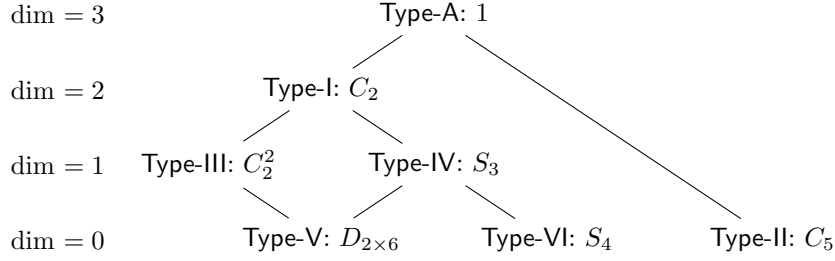
\begin{figure}[ht]
    \centering
    \begin{tikzpicture}
        \node (TypeA) at (0,0) {\TypeA: \(1\)} ;

        \node (TypeII) at (4.5,-3) {\TypeII: \(C_5\)} ;

        \node (TypeI) at (-1.5,-1) {\TypeI: \(C_2\)} ;

        \node (TypeIII) at (-3,-2) {\TypeIII: \(C_{2}^2\)} ;

        \node (TypeIV) at (0,-2) {\TypeIV: \(S_3\)} ;

        \node (TypeV) at (-1.5,-3) {\TypeV: \(D_{2\times 6}\)} ;

        \node (TypeVI) at (1.5,-3) {\TypeVI: \(S_4\)} ;

        \path (TypeII) edge (TypeA)  ;
        \path (TypeI) edge (TypeA)   ;
        \path (TypeIII) edge (TypeI) ;
        \path (TypeIV) edge (TypeI)  ;
        \path (TypeV) edge (TypeIII) ;
        \path (TypeV) edge (TypeIV)  ;
        \path (TypeVI) edge (TypeIV) ;

        \node (Mzero) at (-5,-3) {dim \(= 0\)} ;
        \node (Mone) at (-5,-2) {dim \(= 1\)} ;
        \node (Mtwo) at (-5,-1) {dim \(= 2\)} ;
        \node (Mthree) at (-5,0) {dim \(=3\)} ;
    \end{tikzpicture}
    \caption{Reduced automorphism groups for genus-2 Jacobians.
        Dimensions are of the corresponding loci 
        in the 3-dimensional moduli space of PPASes.
        Lines connect sub- and super-types.
    }
    \label{fig:RA-Jacobian}
\end{figure}

There are
seven possible reduced automorphism groups
for Jacobian surfaces
(provided \(p > 5\); see~\cite{1887/Bolza}).
Figure~\ref{fig:RA-Jacobian}
gives the taxonomy of Jacobian surfaces by reduced automorphism group,
using Bolza's names (``types'') for the classes of Jacobian surfaces
with each of the reduced automorphism groups
(we add \TypeA for the Jacobians with trivial reduced automorphism
group).
We will give normal forms for each type in~\S\ref{sec:atlas}.

We can identify the isomorphism class of a Jacobian
using the Clebsch invariants:
\[
    \classof{\Jac{\XC}}
    \longleftrightarrow
    (A:B:C:D)
    \in \PP(2,4,6,10)(\field)
    \,,
\]
where \(A\), \(B\), \(C\), and \(D\)
are homogeneous polynomials of degree 2, 4, 6, and 10
in the coefficients of the sextic defining~\(\XC\)
(see~\cite[\S1]{1991/Mestre}).
They should be seen as coordinates on the weighted projective space
\(\PP(2,4,6,10)\):
that is, 
\[
    (A:B:C:D) = (\lambda^2 A: \lambda^4 B: \lambda^6 C: \lambda^{10} D)
    \quad
    \text{for all }
    \lambda \not= 0 \in \fieldbar
    \,.
\]
We will not define \((A:B:C:D)\) explicitly here;
in practice, we compute them using
(e.g.) \texttt{ClebschInvariants} in Magma~\cite{2020/Magma}
or \texttt{clebsch\_invariants} in Sage~\cite{2020/SageMath}.

To determine \(\RAut(\Jac{\XC})\) for a given \(\XC\),
we use Bolza's criteria on Clebsch invariants
given in Table~\ref{tab:RAut-from-Clebsch}.
We will need some derived invariants
(see~\cite{1991/Mestre}):
let
\begin{align*}
    A_{11} & = 2C + \frac{1}{3}AB
    \,,
    &
    A_{12} & = \frac{2}{3}(B^2 + AC)
    \,,
    &
    A_{23} & = \frac{1}{2}B\cdot A_{12} + \frac{1}{3}C\cdot A_{11}
    \,,
    \\
    A_{22} & = D
    \,,
    &
    A_{31} & = D
    \,,
    &
    A_{33} & = \frac{1}{2}B\cdot A_{22} + \frac{1}{3}C\cdot A_{12}
    \,,
\end{align*}
and let \(R\) be defined by \(2R^2 = \det (A_{ij})\)
(we will only need to know
whether \(R = 0\)).

\begin{table}[ht]
    \centering
    \begin{tabular}{|c|c|}
        \hline
        Type 
        & Conditions on Clebsch invariants 
        \\
        \hline
        \hline
        \TypeA 
               & \(R \not= 0\), \((A:B:C:D) \not= (0:0:0:1)\)
        \\
        \hline
        \TypeI 
               & \(R = 0\), \(A_{11}A_{22} \not= A_{12}\)
        \\
        \hline
        \TypeII 
                & \((A:B:C:D) = (0:0:0:1)\)
        \\
        \hline
        \TypeIII
                 & \(BA_{11} - 2AA_{12} = -6D\),
                 \(CA_{11} + 2BA_{12} = AD\),
                 \(6C^2 \not= B^3\),
                 \(D \not= 0\)
        \\
        \hline
        \TypeIV 
                 & \(6C^2 = B^3\),
                   \(3D = 2BA_{11}\),
                   \(2AB \not= 15C\),
                   \(D \not= 0\)
        \\
        \hline
        \TypeV   
               & \(6B = A^2\),
                 \(D = 0\),
                 \(A_{11} = 0\),
                 \(A \not= 0\)
        \\
        \hline
        \TypeVI  
                & \((A:B:C:D) = (1:0:0:0)\)
        \\
        \hline
    \end{tabular}
    \caption{%
        Determining the \RAut-type of \(\Jac{\XC}\)
        from its Clebsch invariants.
    }
    \label{tab:RAut-from-Clebsch}
\end{table}

\subsection{Reduced automorphism groups of elliptic products}

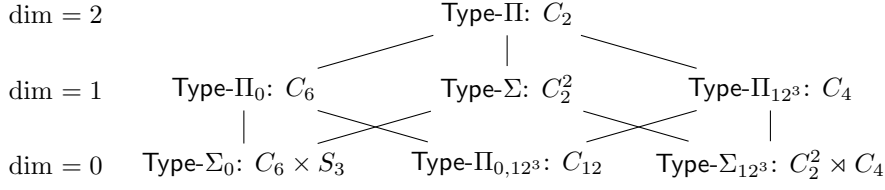
\begin{figure}[ht]
    \centering
    \begin{tikzpicture}

        \node (TypeExE) at (0,0) {\TypeExE: \(C_2\)} ;

        \node (TypeE2) at (0,-1) {\TypeEsquared: \(C_2^2\)} ;

        \node (TypeExE0) at (-3.5,-1) {\TypeExEzero: \(C_6\)} ;

        \node (TypeExE1728) at (3.5,-1) {\TypeExEtwelvecubed: \(C_4\)} ;

        \node (TypeE02) at (-3.5,-2) {\TypeEzerosquared: \(C_6\times S_3\)} ;

        \node (TypeE0xE1728) at (0,-2) {\TypeEzeroxEtwelvecubed: \(C_{12}\)} ;

        \node (TypeE17282) at (3.5,-2) {\TypeEtwelvecubedsquared: \(C_2^2\rtimes C_4\)} ;
                 
        \path (TypeExE0) edge (TypeExE)  ;
        \path (TypeE2) edge (TypeExE)   ;
        \path (TypeExE1728) edge (TypeExE) ;
        \path (TypeE02) edge (TypeExE0)  ;
        \path (TypeE02) edge (TypeE2)  ;
        \path (TypeE17282) edge (TypeExE1728)  ;
        \path (TypeE17282) edge (TypeE2)  ;
        \path (TypeE0xE1728) edge (TypeExE0)  ;
        \path (TypeE0xE1728) edge (TypeExE1728)  ;

        \node (Mzero) at (-6,-2) {dim \(= 0\)} ;
        \node (Mone) at (-6,-1) {dim \(= 1\)} ;
        \node (Mtwo) at (-6,0) {dim \(= 2\)} ;
    \end{tikzpicture}
    \caption{Reduced automorphism groups of elliptic products.
        Dimensions are of the corresponding loci 
        in the 3-dimensional moduli space of PPASes.
        Lines connect sub- and super-types.
    }
    \label{fig:RA-elliptic}
\end{figure}

There are seven 
possible reduced automorphism groups
for elliptic product surfaces~\cite[Prop.~3]{2020/Florit--Smith}.
Figure~\ref{fig:RA-elliptic}
shows the taxonomy of elliptic product
surfaces by reduced automorphism group.
The names (``types'') for the classes of surfaces
are taken from~\cite{2020/Florit--Smith}.

Every elliptic product \(\EC\times\EC'\)
has an involution \(\sigma = [1]\times[-1]\)
in \(\RAut(\EC\times\EC')\).
If \(\EC \cong \EC'\)
then there is also the involution \(\tau\)
exchanging the factors of the product.
The situation is more complicated
if either or both factors are isomorphic to 
one of 
\begin{align*}
    \EC_0: y^2 & = x^3 - 1
    \qquad \text{with}
    &
    \Aut(\EC_{0})
    & = \subgrp{
        \zeta: (x,y) \mapsto (\zeta_3x,-y)
    }
    \cong 
    C_6 
    \intertext{where \(\zeta_3\) is a primitive 3rd root of unity, or}
    \EC_{12^3}: y^2 & = x^3 - x
    \qquad \text{with}
    & 
    \Aut(\EC_{12^3})
    & = \subgrp{
        \iota: (x,y) \mapsto (-x,\sqrt{-1}y)
    }
    \cong 
    C_4
    \,.
\end{align*}
When constructing isogenies, we label the \(2\)-torsion of \(\EC_0\) and
\(\EC_{12^3}\) as follows:
\begin{align*}
    \EC_{0}[2] 
    & = 
    \{
        0
        ,
        P_1 = (1,0)
        ,
        P_2 = (\zeta_3,0)
        ,
        P_3 = (\zeta_3^2,0)
    \}
    \,,
    \\
    \EC_{12^3}[2] 
    & = 
    \{
        0
        ,
        P_1 = (1,0)
        ,
        P_2 = (-1,0)
        ,
        P_3 = (0,0)
    \}
    \,.
\end{align*}


When navigating isogeny graphs,
we can identify the isomorphism class
of an elliptic product
using the pair of \(j\)-invariants of the factors:
\[
    \classof{\EC_1\times\EC_2}
    \longleftrightarrow
    \{j(\EC_1),j(\EC_2)\}
    \,.
\]
To determine \(\RAut(\EC_1\times\EC_2)\),
we can use the criteria on \(j\)-invariants
given in Table~\ref{tab:RAut-from-j-invariants}.

\begin{table}[ht]
    \centering
    \begin{tabular}{|c|c||c|c|}
        \hline
        Type 
        & Conditions 
        &
        Type 
        & Conditions
        \\
        \hline
        \hline
        \TypeExE
        & \(\{j(\EC_1),j(\EC_2)\}\cap\{0,1728\} = \emptyset\)
        &
        \multirow{2}{*}{\TypeEsquared}
        & \(j(\EC_1) = j(\EC_2)\),
        \\
        \cline{1-2}
        \TypeExEzero 
        & \(j(\EC_1) = 0\) or \(j(\EC_2) = 0\) 
        &
        {}
        & \(j(\EC_i) \not\in \{0,1728\}\)
        \\
        \hline
        \TypeExEtwelvecubed 
        & \(j(\EC_1) = 1728\) or \(j(\EC_2) = 1728\)
        &
        \TypeEzerosquared
        & \(j(\EC_1) = j(\EC_2) = 0\)
        \\
        \hline
        \TypeEzeroxEtwelvecubed 
        & \(\{j(\EC_1),j(\EC_2)\} = \{0,1728\}\)
        &
        \TypeEtwelvecubedsquared 
        & \(j(\EC_1) = j(\EC_2) = 1728\)
        \\
        \hline
    \end{tabular}
    \caption{%
        Determining the \RAut-type of an elliptic product
        \(\EC_1\times\EC_2\).
    }
    \label{tab:RAut-from-j-invariants}
\end{table}


\subsection{Superspecial vertices}
\label{sec:superspecial-count}

Ibukiyama, Katsura, and Oort have computed the precise number of superspecial
genus-2 Jacobians (up to isomorphism) of each reduced automorphism
type~\cite[Theorem 3.3]{1986/Ibukiyama--Katsura--Oort}.
We reproduce their results for \(p > 5\),
completing them with the number of superspecial elliptic products
of each automorphism type
(which can be easily derived from the well-known formula
for the number of supersingular elliptic curves over~\(\FF_{p^2}\))
in Table~\ref{tab:number-of-superspecial-vertices}.

\begin{definition}
    \label{def:epsilons}
    For each prime \(p > 5\), we define the following quantities:
    \begin{itemize}
        \item
            \(\epsilon_{1,p} = 1\) if \(p \equiv 3\pmod{4}\), 0 otherwise;
        \item
            \(\epsilon_{2,p} = 1\) if \(p \equiv 5,7\pmod{8}\), 0 otherwise;
        \item
            \(\epsilon_{3,p} = 1\) if \(p \equiv 2\pmod{3}\), 0 otherwise;
        \item
            \(\epsilon_{5,p} = 1\) if \(p \equiv 4 \pmod{5}\), 0 otherwise;
        \item
            \(N_p = (p-1)/12 - \epsilon_{1,p}/2 - \epsilon_{3,p}/3\).
    \end{itemize}
    Note that \(N_p\), \(\epsilon_{1,p}\), and \(\epsilon_{3,p}\)
    count the isomorphism classes of supersingular elliptic curves
    over \(\FF_{p^2}\)
    with reduced automorphism group of order \(1\), \(2\), and \(3\),
    respectively.
\end{definition}

\begin{table}[ht]
    \centering
    \begin{tabular}{|r|l||r|l|}
        \hline
        Type & Vertices in \(\RichelotIG\)
        &
        Type & Vertices in \(\RichelotIG\)
        \\
        \hline
        \hline
        \multirow{2}{*}{\TypeI}
               & \( \frac{1}{48}(p-1)(p-17) \)
        &       
        \TypeExE & \( \frac{1}{2}N_p(N_p - 1) \)
        \\
        \cline{3-4}
        {}
               & \qquad \quad \(
                   + \frac{1}{4}\epsilon_{1,p}
                   + \epsilon_{2,p}
                   + \epsilon_{3,p}
                 \)
        &
        \TypeExEzero & \( \epsilon_{3,p}N_p \)
        \\
        \hline
        \TypeII & \(\epsilon_{5,p}\)
        &
        \TypeExEtwelvecubed & \( \epsilon_{1,p}N_p \)
        \\
        \hline
        \TypeIII & \( 
                      \frac{3}{2}N_p 
                      + \frac{1}{2}\epsilon_{1,p} 
                      - \frac{1}{2}\epsilon_{2,p} 
                      - \frac{1}{2}\epsilon_{3,p}
                   \)
        &
        \TypeEzeroxEtwelvecubed & \(\epsilon_{1,p}\cdot\epsilon_{3,p}\)
        \\
        \hline
        \TypeIV & \(2N_p + \epsilon_{1,p} - \epsilon_{2,p}\)
        &
        \TypeEsquared & \(N_p\) 
        \\
        \hline
        \TypeV & \(\epsilon_{3,p}\)
        &
        \TypeEzerosquared & \(\epsilon_{3,p}\)
        \\
        \hline
        \TypeVI & \(\epsilon_{2,p}\)
        &
        \TypeEtwelvecubedsquared & \(\epsilon_{1,p}\)
        \\
        \hline
        \TypeA & \multicolumn{3}{l|}{\(
                   \frac{1}{2880}(p-1)(p^2-35p+346)
                   - \frac{1}{16}\epsilon_{1,p}
                   - \frac{1}{4}\epsilon_{2,p}
                   - \frac{2}{9}\epsilon_{3,p}
                   - \frac{1}{5}\epsilon_{5,p}
                 \)}
        \\
        \hline
    \end{tabular}
    \caption{%
        The number of superspecial vertices 
        of each \(\RAut\)-type.
    }
    \label{tab:number-of-superspecial-vertices}
\end{table}

If the reader chooses suitable values of \(p\)
and computes \(\RichelotIG\),
then they will find graphs built from overlapping
copies of the neighbourhoods described in~\S\ref{sec:atlas}.
We will see that \(\RichelotIG\)
is much more complicated
than the elliptic \(2\)-isogeny graph.

\section{
    An atlas of the Richelot isogeny graph
}
\label{sec:atlas}

We are now ready to compute the neighbourhoods of each type of vertex
and edge in the Richelot isogeny graph.
We begin with general (\TypeA) vertices,
before considering each type with an involution,
in order of increasing speciality,
and ending with \TypeII (which has no involution).

\subsection{The algorithm}
We compute each vertex neighbourhood in the same way:
\begin{enumerate}
    \item
        Take the generic curve or product for the \RAut-type.
        We use Bolza's 
        normal forms for the curves with special reduced automorphism
        groups from Bolza~\cite{1887/Bolza}, reparametrizing 
        to force full rational \(2\)-torsion in the Jacobians.
    \item
        Enumerate the \((2,2)\)-isogeny kernels.
    \item
        Compute the action of the reduced automorphism group.
    \item
        For each orbit,
        choose a representative kernel,
        compute the codomain
        using the formul\ae{} of~\S\ref{sec:Richelot}
        and~\S\ref{sec:elliptic-product-isogenies},
        and identify the \RAut-type of the codomain
        using the criteria of Tables~\ref{tab:RAut-from-Clebsch}
        and~\ref{tab:RAut-from-j-invariants}.
        The orbit sizes give edge weights.
\end{enumerate}
For subsequent isogenies,
we repeat Steps (2), (3), and (4)
from the current vertex.

\subsection{Diagram notation}

In all of our diagrams, 
\textbf{solid vertices} have definite types,
and \textbf{solid edges} have definite weights.
The \textbf{dotted vertices} have an indicative type,
but may change type under specialization,
acquiring more automorphisms,
with the weight of \textbf{dotted edges} increasing proportionally
according to Eq.~\eqref{eq:ratio-principle}.
For example: in Figure~\ref{fig:Type-A},
if one of the dotted neighbours specializes to a \TypeI
vertex, then the returning dotted arrow will become a weight-2 arrow.
All edges from solid vertices are shown;
some edges from dotted vertices, especially to vertices outside the diagram,
are omitted for clarity.

\subsection{General vertices and edges}

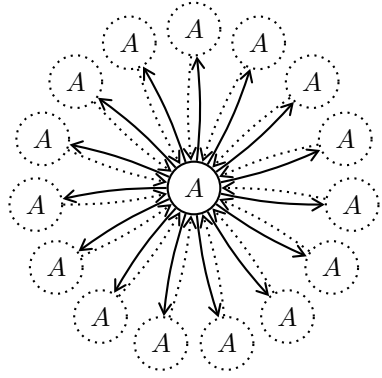
\begin{figure}[ht]
    \centering
    \begin{tikzpicture}[
            >={angle 60},
            thick,
            vertex/.style = {circle, draw, inner sep=0.5mm, minimum size=7mm}
        ]
        \node (s) [
            regular polygon,
            regular polygon sides=15,
            minimum size=42mm,
            above]
        at (0,0) {};
        \node (c) [vertex] at (s.center) {$A$};
        \foreach \i in {1,...,15}{
            \node (s\i) [vertex, dotted] at (s.corner \i) {$A$} ;
            \draw[->] (c) edge[bend right=6] (s\i) ;
            \draw[->] (s\i) edge[bend right=6, dotted] (c) ;
        }

    \end{tikzpicture}
    \caption{The neighbourhood of a \protect{\TypeA} vertex.}
    \label{fig:Type-A}
\end{figure}

Figure~\ref{fig:Type-A} shows the neighbourhood of a \TypeA vertex:
there are weight-1 edges to fifteen neighbouring vertices, generally all \TypeA,
and a weight-1 dual edge returning from each of them.

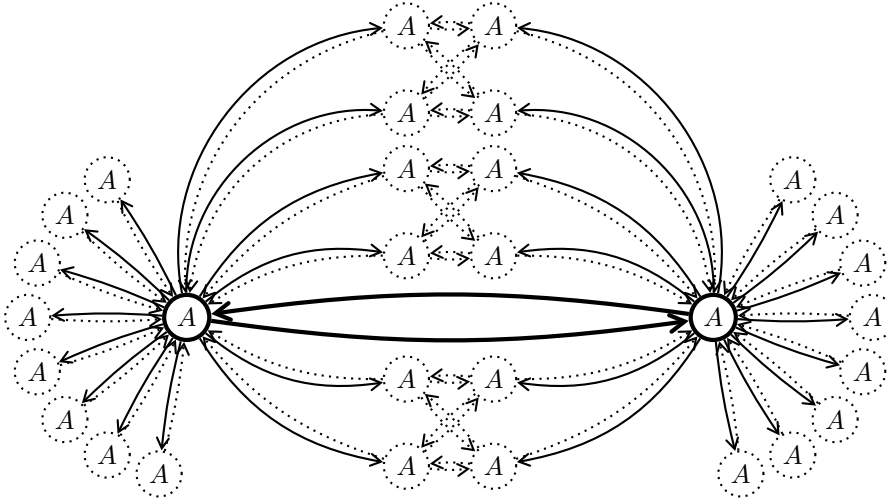
\begin{figure}[ht]
    \centering
    \begin{tikzpicture}[
            >={angle 60},
            thick,
            vertex/.style = {circle, draw, fill=white, inner sep=0.5mm, minimum size=6mm}
        ]
        \node (domain) [vertex, ultra thick] at (0,0) {$A$} ;
        \node (codomain) [vertex, ultra thick] at (7,0) {$A$} ;
        \draw[->, ultra thick] (codomain) edge[bend right=8] (domain) ;
        \draw[->, ultra thick] (domain) edge[bend right=8] (codomain) ;

        \node (Xtop) [
            regular polygon,
            regular polygon sides=4,
            minimum size=16mm,
        ] at (3.5, 3.3) {};
        \node (Xtopcod1) [vertex, dotted] at (Xtop.corner 1) {$A$} ;
        \node (Xtopdom1) [vertex, dotted] at (Xtop.corner 2) {$A$} ;
        \node (Xtopdom2) [vertex, dotted] at (Xtop.corner 3) {$A$} ;
        \node (Xtopcod2) [vertex, dotted] at (Xtop.corner 4) {$A$} ;
        \foreach \i in {1,2} {
            \draw[->] (domain) edge[bend right=-48] (Xtopdom\i) ;
            \draw[->] (Xtopdom\i) edge[bend right=40, dotted] (domain) ;
            \draw[->] (codomain) edge[bend right=48] (Xtopcod\i) ;
            \draw[->] (Xtopcod\i) edge[bend right=-40, dotted] (codomain) ;
            \foreach \j in {1,2} {
                \draw[->] (Xtopdom\i) edge[bend right=6, dotted] (Xtopcod\j) ;
                \draw[->] (Xtopcod\j) edge[bend right=6, dotted] (Xtopdom\i) ;
            }
        }

        \node (Xmid) [
            regular polygon,
            regular polygon sides=4,
            minimum size=16mm,
        ] at (3.5, 1.4) {};
        \node (Xmidcod1) [vertex, dotted] at (Xmid.corner 1) {$A$} ;
        \node (Xmiddom1) [vertex, dotted] at (Xmid.corner 2) {$A$} ;
        \node (Xmiddom2) [vertex, dotted] at (Xmid.corner 3) {$A$} ;
        \node (Xmidcod2) [vertex, dotted] at (Xmid.corner 4) {$A$} ;
        \foreach \i in {1,2} {
            \draw[->] (domain) edge[bend right=-24] (Xmiddom\i) ;
            \draw[->] (Xmiddom\i) edge[bend right=16, dotted] (domain) ;
            \draw[->] (codomain) edge[bend right=24] (Xmidcod\i) ;
            \draw[->] (Xmidcod\i) edge[bend right=-16, dotted] (codomain) ;
            \foreach \j in {1,2} {
                \draw[->] (Xmiddom\i) edge[bend right=6, dotted] (Xmidcod\j) ;
                \draw[->] (Xmidcod\j) edge[bend right=6, dotted] (Xmiddom\i) ;
            }
        }

        \node (Xbot) [
            regular polygon,
            regular polygon sides=4,
            minimum size=16mm,
        ] at (3.5, -1.4) {};
        \node (Xbotcod1) [vertex, dotted] at (Xbot.corner 1) {$A$} ;
        \node (Xbotdom1) [vertex, dotted] at (Xbot.corner 2) {$A$} ;
        \node (Xbotdom2) [vertex, dotted] at (Xbot.corner 3) {$A$} ;
        \node (Xbotcod2) [vertex, dotted] at (Xbot.corner 4) {$A$} ;
        \foreach \i in {1,2} {
            \draw[->] (domain) edge[bend right=24] (Xbotdom\i) ;
            \draw[->] (Xbotdom\i) edge[bend right=-16, dotted] (domain) ;
            \draw[->] (codomain) edge[bend right=-24] (Xbotcod\i) ;
            \draw[->] (Xbotcod\i) edge[bend right=16, dotted] (codomain) ;
            \foreach \j in {1,2} {
                \draw[->] (Xbotdom\i) edge[bend right=6, dotted] (Xbotcod\j) ;
                \draw[->] (Xbotcod\j) edge[bend right=6, dotted] (Xbotdom\i) ;
            }
        }

        \node (domainin) [
            regular polygon,
            regular polygon sides=18,
            minimum size=42mm,
        ] at (domain) {};
        \foreach \i in {3,4,5,6,7,8,9,10} {
            \node (domA\i) [vertex, dotted] at (domainin.corner \i) {$A$} ;
            \draw[->] (domain) edge[bend right=4] (domA\i) ;
            \draw[->] (domA\i) edge[bend right=4, dotted] (domain) ;
        }

        \node (codomainout) [
            regular polygon,
            regular polygon sides=18,
            minimum size=42mm,
        ] at (codomain) {};
        \foreach \i in {11,12,13,14,15,16,17,18} {
            \node (codA\i) [vertex, dotted] at (codomainout.corner \i) {$A$} ;
            \draw[->] (codomain) edge[bend right=4] (codA\i) ;
            \draw[->] (codA\i) edge[bend right=4, dotted] (codomain) ;
        }
    \end{tikzpicture}
    \caption{The neighbourhood of a general edge and its dual.}
    \label{fig:general-edge}
\end{figure}

The Richelot isogeny graph is 15-regular (counting weights),
and it is tempting to imagine that locally, the graph looks 
like an assembly of copies of the star in Figure~\ref{fig:Type-A},
with each outer vertex becoming the centre of its own star.
However,
the reality is more complicated.
If we look at a pair of neighbouring \TypeA vertices,
then six of the neighbours of one are connected to neighbours of the
other.
Figure~\ref{fig:general-edge} shows this configuration.

The interconnections in Figure~\ref{fig:general-edge}
are explained as follows.
For each \((2,2)\)-isogeny \(\phi: \AV_1 \to \AV_{2}\),
there are \emph{twelve} \((4,4,2,2)\)-isogenies (each a composition of
three \((2,2)\)-isogenies) from \(\AV_{1}\) to \(\AV_{2}\);
composing any of these with \(\dualof{\phi}\)
defines a cycle of length 4 in the graph,
which is isomorphic to multiplication-by-4 on \(\AV_1\).
These cycles of length 4 are the ``small cycles'' 
exploited
by Flynn and Ti in~\cite[\S2.3]{2019/Flynn--Ti}.
In contrast, composing a central isogeny
with one of the eight isogenies from the far left
or the eight from the far right 
of Figure~\ref{fig:general-edge}
yields a \((4,4)\)-isogeny,
and composing with one of each yields an \((8,8)\)-isogeny.
In the terminology of~\cite{2020/Castryck--Decru--Smith},
the isogenies at the far left and far right are ``good'' extensions of
the central pair, while those forming the adjacent edges of squares
are ``bad'' extensions of each other.

This pattern is replicated throughout the Richelot isogeny graph:
each edge is common to twelve of these 4-cycles (counting weights as
multiplicities).

\subsection{General elliptic products:
\texorpdfstring{\TypeExE}{Type-Pi} vertices}

\begin{figure}[ht]
    \centering
    \begin{tabular}{c@{\qquad}c}
        \begin{tikzpicture}[
                >={angle 60},
                thick,
                vertex/.style = {circle, draw, fill=white, inner sep=0.5mm, minimum size=7mm}
            ]
            \node (s) [
                regular polygon,
                regular polygon sides=15,
                minimum size=42mm,
                above
            ] at (0,0) {};
            \node (c) [vertex] at (s.center) {\TypeExElabel};
            \foreach \i in {1,2,5,6,7,10,11,12,15}{
                \node (s\i) [vertex, dotted] at (s.corner \i)
                {\TypeExElabel} ;
                \draw[->] (c) edge[bend right=8] (s\i) ;
                \draw[->] (s\i) edge[bend right=8, dotted] (c) ;
            }
            \foreach \i in {3,4,8,9,13,14}{
                \node (s\i) [vertex, dotted] at (s.corner \i) {$I$} ;
                \draw[->] (c) edge[bend right=8] (s\i) ;
                \draw[->] (s\i) edge[bend right=8, dotted] (c) ;
            }
        \end{tikzpicture}
        &
        \begin{tikzpicture}[ 
                >={angle 60},
                thick,
                vertex/.style = {circle, draw, fill=white, inner sep=0.5mm, minimum size=7mm},
                wt/.style = {fill=white, anchor=center, pos=0.5, minimum size=1mm, inner sep=1pt}
            ]
            \node (s) [
                regular polygon,
                regular polygon sides=11,
                minimum size=42mm,
                above
            ] at (0,0) {};
            \node (c) [vertex] at (s.center) {$I$};
            \node (s1) [vertex,dotted] at (s.corner 1) {\TypeExElabel};
            \draw[->] (c) edge[bend right=8] (s1) ;
            \draw[->] (s1) edge[bend right=8, dotted] (c) ;

            \foreach \i in {2,3,4,9,10,11}{
                \node (s\i) [vertex,dotted] at (s.corner \i) {$I$}; 
                \draw[->] (c) edge[bend right=8] (s\i) ;
                \draw[->] (s\i) edge[bend right=8, dotted] (c) ;
            }
            \foreach \i in {5,6,7,8}{
                \node (s\i) [vertex,dotted] at (s.corner \i) {$A$}; 
                \draw[->] (c) edge[bend right=8] node[wt]{2}(s\i) ;
                \draw[->] (s\i) edge[bend right=8, dotted] (c) ;
            }
        \end{tikzpicture}
    \end{tabular}
    \caption{Neighbourhoods of the general \TypeExE and \TypeI vertices.}
    \label{fig:TypeExE-TypeI}
\end{figure}
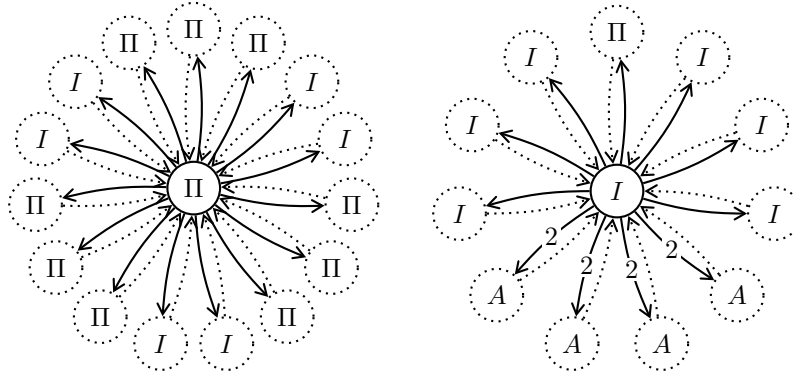

The general \TypeExE vertex
is an elliptic product vertex \(\classof{\EC\times\EC'}\)
where \(\EC' \not\cong \EC\),
and neither \(\EC\) nor \(\EC'\)
has special automorphisms.
In this case
\(\RAut(\EC\times\EC') = \subgrp{\sigma} \cong C_2\),
which fixes every \((2,2)\)-isogeny kernel,
so we have a subgroup isomorphic to \(C_2\) in the reduced automorphism group
of every \((2,2)\)-isogeny codomain.
The nine elliptic product neighbours are generally \TypeExE;
the six Jacobian neighbours are generally \TypeI,
the most general type with a reduced involution.
The situation is illustrated at the left of Figure~\ref{fig:TypeExE-TypeI}.

\begin{remark}
    Looking at Figure~\ref{fig:TypeExE-TypeI},
    we see that
    \TypeExE vertices cannot have \TypeA or \TypeII neighbours:
    any walk in the graph from a \TypeA vertex to an elliptic product
    must have already passed through a vertex with an involution
    in its reduced automorphism group.
    We will see below that the same applies
    to any elliptic product or square vertex,
    as well as to \TypeIV, \TypeV, and \TypeVI vertices.
\end{remark}

\subsection{\texorpdfstring{\TypeI}{Type-I} vertices}
\label{sec:TypeI}

The generic \TypeI vertex is \(\classof{\Jac{\XC_{I}}}\),
where \(\XC_{I}\) 
is defined by
\[
    \XC_{I}: y^2 = F_{I}(x) := (x^2 - 1)(x^2 - s^2)(x^2 - t^2)
\]
with parameters \(s\) and \(t\).
Any Jacobian \(\AV_{0}\) with \(C_2 \subseteq \RAut(\AV_{0})\)
(that is, \TypeI, \TypeIII, \TypeIV, \TypeV, or \TypeVI)
is isomorphic to the Jacobian of \(\Jac{\XC_I}\)
for some \((s,t)\)
such that
\(
    st(s^2-1)(t^2-1)(s^2-t^2)
    \not=
    0
\).

There are maximal \(2\)-Weil isotropic subgroups 
\(K_{1},\ldots,K_{15}\) of \(\Jac{\XC_I}[2]\);
each is the kernel of a \((2,2)\)-isogeny \(\Jac{\XC_I} = \AV_{0} \to \AV_i = \AV_{0}/K_i\).  
The kernels \(K_{i}\) correspond to the following quadratic splittings.
First:
\begin{align*}
    K_{1} \leftrightarrow \{ x^2-1, x^2-s^2, x^2-t^2 \}
    \,.
\end{align*}
These three quadratics are linearly dependent,
so \(\AV_{1} \cong \EC\times\EC'\)
with factors
\(\EC: y^2 = (x-1)(x-s^2)(x-t^2)\)
and
\(\EC': y^2 = (x-1)(x-1/s^2)(x-1/t^2)\).

Six of the kernels share a nontrivial element with \(K_{1}\),
namely
\begin{align*}
    K_{2} &\leftrightarrow \{ x^2 - 1, x^2 \pm (s + t)x + st \}
    \,,
    &
    K_{3} &\leftrightarrow \{ x^2 - 1, x^2 \pm (s - t)x - st \}
    \,,
    \\
    K_{4} &\leftrightarrow \{ x^2 - s^2, x^2 \pm (t+1)x + t \}
    \,,
    &
    K_{5} &\leftrightarrow \{ x^2 - s^2, x^2 \pm (t-1)x - t \}
    \,,
    \\
    K_{6} &\leftrightarrow \{ x^2 - t^2, x^2 \pm (s+1)x + s \}
    \,,
    &
    K_{7} &\leftrightarrow \{ x^2 - t^2, x^2 \pm (s-1)x - s \}
    \,.
\end{align*}
The last eight kernels do not share any nontrivial
elements with \(K_{1}\), namely
\begin{align*}
    K_{8} &\leftrightarrow
        \{
            x^2 + (s - 1)x - s,
            x^2 - (t - 1)x - t,
            x^2 - (s - t)x - st
        \}
    \,,
    \\
    K_{9} &\leftrightarrow
        \{
            x^2 - (s - 1)x - s,
            x^2 + (t - 1)x - t,
            x^2 + (s - t)x - st
        \}
    \,,
    \\
    K_{10} &\leftrightarrow
        \{
            x^2 - (s - 1)x - s,
            x^2 - (t + 1)x + t,
            x^2 + (s + t)x + st
        \}
    \,,
    \\
    K_{11} &\leftrightarrow
        \{
            x^2 + (s - 1)x - s,
            x^2 + (t + 1)x + t,
            x^2 - (s + t)x + st
        \}
    \,,
    \\
    K_{12} &\leftrightarrow
        \{
            x^2 + (s + 1)x + s,
            x^2 + (t - 1)x - t,
            x^2 - (s + t)x + st
        \}
    \,,
    \\
    K_{13} &\leftrightarrow
        \{
            x^2 - (s + 1)x + s,
            x^2 - (t - 1)x - t,
            x^2 + (s + t)x + st
        \}
    \,,
    \\
    K_{14} &\leftrightarrow
        \{
            x^2 + (s + 1)x + s,
            x^2 - (t + 1)x + t,
            x^2 - (s - t)x - st
        \}
    \,,
    \\
    K_{15} &\leftrightarrow
        \{
            x^2 - (s + 1)x + s,
            x^2 + (t + 1)x + t,
            x^2 + (s - t)x - st
        \}
    \,.
\end{align*}

The reduced automorphism group is 
\(\RAut(\XC_{I}) = \subgrp{\sigma} \cong C_2\),
where \(\sigma\) acts as
\[
    \sigma_*: x \longleftrightarrow -x
\]
on \(x\)-coordinates, and 
on (the indices of) the set of kernels \(\{K_{1},\ldots,K_{15}\}\)
via
\begin{align*}
    \sigma_*
    =
    (1)(2)(3)(4)(5)(6)(7)(8,9)(10,11)(12,13)(14,15)
    \,.
\end{align*}
The orbits of the kernel subgroups under \(\sigma\)
and the types of the corresponding neighbours
are listed in Table~\ref{tab:TypeI}.
The situation is illustrated on the right of
Figure~\ref{fig:TypeExE-TypeI}.

\begin{table}[ht]
    \centering
    \begin{tabular}{c|c|c||c|c|c}
        Kernel orbit
        & Stabilizer
        & Codomain
        & Kernel orbit
        & Stabilizer
        & Codomain
        \\
        \hline
        \(\{K_{1}\}\) & \(\subgrp{\sigma}\) & \TypeExE
        &
        \(\{K_{7}\}\) & \(\subgrp{\sigma}\) & \TypeI
        \\
        \(\{K_{2}\}\) & \(\subgrp{\sigma}\) & \TypeI
        &
        \(\{K_{8,9}\}\) & \(1\) & \TypeA
        \\
        \(\{K_{3}\}\) & \(\subgrp{\sigma}\) & \TypeI
        &
        \(\{K_{10,11}\}\) & \(1\) & \TypeA
        \\
        \(\{K_{4}\}\) & \(\subgrp{\sigma}\) & \TypeI
        &
        \(\{K_{12,13}\}\) & \(1\) & \TypeA
        \\
        \(\{K_{5}\}\) & \(\subgrp{\sigma}\) & \TypeI
        &
        \(\{K_{14,15}\}\) & \(1\) & \TypeA
        \\
        \(\{K_{6}\}\) & \(\subgrp{\sigma}\) & \TypeI
    \end{tabular}
    \caption{%
        Edge data for the generic \TypeI vertex.
    }
    \label{tab:TypeI}
\end{table}
 
\begin{figure}[ht]
    \centering
    \begin{tikzpicture}[
            >={angle 60},
            thick,
            vertex/.style = {circle, draw, fill=white, inner sep=0.5mm, minimum size=6mm},
            wt/.style = {fill=white, anchor=center, pos=0.5, minimum size=1mm, inner sep=1pt}
        ]
        \node (domain) [vertex, ultra thick] at (0,0) {$I$} ;
        \node (codomain) [vertex, ultra thick] at (7,0) {\TypeExElabel} ;
        \draw[->, ultra thick] (codomain) edge[bend right=8] (domain) ;
        \draw[->, ultra thick] (domain) edge[bend right=8] (codomain) ;

        \node (Xtop) [
            regular polygon,
            regular polygon sides=4,
            minimum size=16mm,
        ] at (3.5, 3.5) {};
        \node (Xtopcod1) [vertex, dotted] at (Xtop.corner 1) {\TypeExElabel} ;
        \node (Xtopdom1) [vertex, dotted] at (Xtop.corner 2) {$I$} ;
        \node (Xtopdom2) [vertex, dotted] at (Xtop.corner 3) {$I$} ;
        \node (Xtopcod2) [vertex, dotted] at (Xtop.corner 4) {$I$} ;
        \foreach \i in {1,2} {
            \draw[->] (domain) edge[bend right=-48] (Xtopdom\i) ;
            \draw[->] (Xtopdom\i) edge[bend right=40, dotted] (domain) ;
            \draw[->] (codomain) edge[bend right=48] (Xtopcod\i) ;
            \draw[->] (Xtopcod\i) edge[bend right=-40, dotted] (codomain) ;
            \foreach \j in {1,2} {
                \draw[->] (Xtopdom\i) edge[bend right=6, dotted] (Xtopcod\j) ;
                \draw[->] (Xtopcod\j) edge[bend right=6, dotted] (Xtopdom\i) ;
            }
        }

        \node (Xmid) [
            regular polygon,
            regular polygon sides=4,
            minimum size=16mm,
        ] at (3.5, 1.5) {};
        \node (Xmidcod1) [vertex, dotted] at (Xmid.corner 1) {\TypeExElabel} ;
        \node (Xmiddom1) [vertex, dotted] at (Xmid.corner 2) {$I$} ;
        \node (Xmiddom2) [vertex, dotted] at (Xmid.corner 3) {$I$} ;
        \node (Xmidcod2) [vertex, dotted] at (Xmid.corner 4) {$I$} ;
        \foreach \i in {1,2} {
            \draw[->] (domain) edge[bend right=-24] (Xmiddom\i) ;
            \draw[->] (Xmiddom\i) edge[bend right=16, dotted] (domain) ;
            \draw[->] (codomain) edge[bend right=24] (Xmidcod\i) ;
            \draw[->] (Xmidcod\i) edge[bend right=-16, dotted] (codomain) ;
            \foreach \j in {1,2} {
                \draw[->] (Xmiddom\i) edge[bend right=6, dotted] (Xmidcod\j) ;
                \draw[->] (Xmidcod\j) edge[bend right=6, dotted] (Xmiddom\i) ;
            }
        }

        \node (Xbot) [
            regular polygon,
            regular polygon sides=4,
            minimum size=16mm,
        ] at (3.5, -1.5) {};
        \node (Xbotcod1) [vertex, dotted] at (Xbot.corner 1) {\TypeExElabel} ;
        \node (Xbotdom1) [vertex, dotted] at (Xbot.corner 2) {$I$} ;
        \node (Xbotdom2) [vertex, dotted] at (Xbot.corner 3) {$I$} ;
        \node (Xbotcod2) [vertex, dotted] at (Xbot.corner 4) {$I$} ;
        \foreach \i in {1,2} {
            \draw[->] (domain) edge[bend right=24] (Xbotdom\i) ;
            \draw[->] (Xbotdom\i) edge[bend right=-16, dotted] (domain) ;
            \draw[->] (codomain) edge[bend right=-24] (Xbotcod\i) ;
            \draw[->] (Xbotcod\i) edge[bend right=16, dotted] (codomain) ;
            \foreach \j in {1,2} {
                \draw[->] (Xbotdom\i) edge[bend right=6, dotted] (Xbotcod\j) ;
                \draw[->] (Xbotcod\j) edge[bend right=6, dotted] (Xbotdom\i) ;
            }
        }

        \node (domainin) [
            regular polygon,
            regular polygon sides=13,
            minimum size=42mm,
        ] at (domain) {};
        \foreach \i in {3,4,5,6} {
            \node (domA\i) [vertex, dotted] at (domainin.corner \i) {$A$} ;
            \draw[->] (domain) edge[bend right=8] node[wt]{2} (domA\i) ;
            \draw[->] (domA\i) edge[bend right=8, dotted] (domain) ;
        }

        \node (codomainout) [
            regular polygon,
            regular polygon sides=20,
            minimum size=42mm,
        ] at (codomain) {};
        \foreach \i in {13,20} {
            \node (codA\i) [vertex, dotted] at (codomainout.corner \i) {$I$} ;
            \draw[->] (codomain) edge[bend right=4] (codA\i) ;
            \draw[->] (codA\i) edge[bend right=4, dotted] (codomain) ;
        }
        \foreach \i in {14,15,16,17,18,19} {
            \node (codA\i) [vertex, dotted] at (codomainout.corner \i) {\TypeExElabel} ;
            \draw[->] (codomain) edge[bend right=4] (codA\i) ;
            \draw[->] (codA\i) edge[bend right=4, dotted] (codomain) ;
        }

    \end{tikzpicture}
    \caption{The neighbourhood of a 
    general \TypeI vertex and its \TypeExE neighbour.}
    \label{fig:TypeI-TypeExE-connect}
\end{figure}
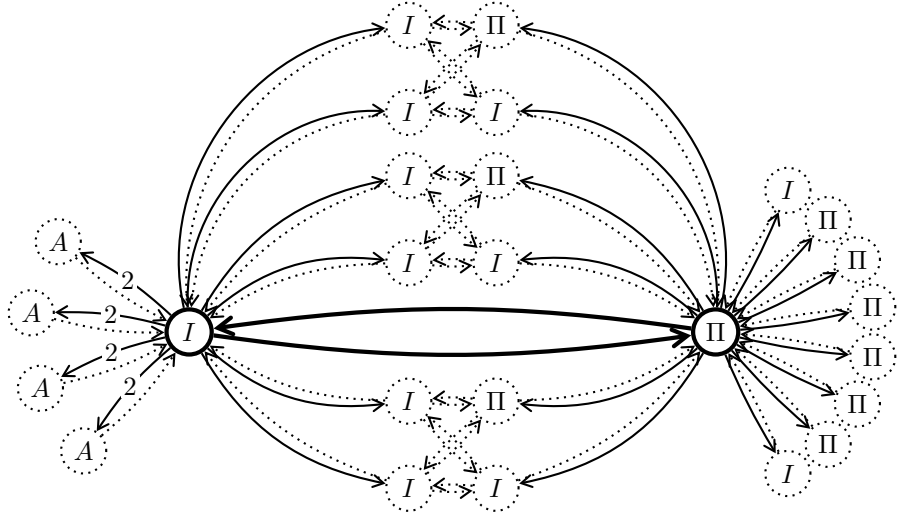

Computing one isogeny step beyond each \TypeI neighbour
of \(\classof{\Jac{\XC_{I}}}\),
we find six neighbours of \(\classof{\EC\times\EC'}\);
thus we complete Figure~\ref{fig:TypeI-TypeExE-connect},
which shows the neighbourhood of the edge \(\classof{\phi_{1}}\)
and its dual, \(\classof{\dualof{\phi_{1}}} = \classof{\phi_{Id}}\).
This should be compared with Figure~\ref{fig:general-edge}.
Note that \(\phi_{i}\circ\dualof{\phi_{1}}\) 
is a \((4,2,2)\)- resp.~\((4,4)\)-isogeny 
for \(2 \le i \le 7\) resp. \(8 \le i \le 15\).

\subsection{General elliptic squares:
\texorpdfstring{\TypeEsquared}{Type-Sigma} vertices}

\begin{figure}[ht]
    \centering
    \begin{tabular}{c@{\qquad}c}
        \begin{tikzpicture}[
                >={angle 60},
                thick,
                vertex/.style = {circle, draw, fill=white, inner sep=0.5mm, minimum size=7mm},
                wt/.style = {fill=white, anchor=center, pos=0.5, minimum size=1mm, inner sep=1pt}
            ]
            \node (s) [
                regular polygon,
                regular polygon sides=11,
                minimum size=42mm,
                above] at (0,0) {};
            \node (c) [vertex] at (s.center) {\TypeEsquaredlabel};
            \foreach \i in {1,4,7}{
                \node (s\i) [vertex, dotted] at (s.corner \i) {\TypeExElabel} ;
                \draw[->] (c) edge[bend right=12] node[wt]{2} (s\i) ;
                \draw[->] (s\i) edge[bend right=12, dotted] (c) ;
            }
            \foreach \i in {2,5,8}{
                \node (s\i) [vertex, dotted] at (s.corner \i) {\TypeEsquaredlabel} ;
                \draw[->] (c) edge[bend right=12] (s\i) ;
                \draw[->] (s\i) edge[bend right=12, dotted] (c) ;
            }
            \foreach \i in {3,6,9}{
                \node (s\i) [vertex, dotted] at (s.corner \i) {$III$} ;
                \draw[->] (c) edge[bend right=12] (s\i) ;
                \draw[->] (s\i) edge[bend right=12, dotted] (c) ;
            }
            \draw[->] (s2) edge[bend right=48, dotted] (s3) ;
            \draw[->] (s3) edge[bend right=-24, dotted] (s2) ;
            \draw[->] (s5) edge[bend right=48, dotted] (s6) ;
            \draw[->] (s6) edge[bend right=-24, dotted] (s5) ;
            \draw[->] (s8) edge[bend right=48, dotted] (s9) ;
            \draw[->] (s9) edge[bend right=-24, dotted] (s8) ;

            \node (s10) [vertex, dotted] at (s.corner 10) {$I$} ;
            \draw[->] (c) edge[bend right=12] node[wt]{2} (s10) ;
            \draw[->] (s10) edge[bend right=12, dotted] (c) ;

            \draw[->] (c) edge[out=45,in=65,loop, looseness=30] (c) ;

        \end{tikzpicture}
        &
        \begin{tikzpicture}[
                >={angle 60},
                thick,
                vertex/.style = {circle, draw, fill=white, inner sep=0.5mm, minimum size=7mm},
                wt/.style = {fill=white, anchor=center, pos=0.5, minimum size=1mm, inner sep=1pt}
            ]
            \node (s) [
                regular polygon,
                regular polygon sides=8,
                minimum size=42mm,
                above
            ] at (0,0) {};
            \node (c) [vertex] at (s.center) {$III$};
            \foreach \i in {1,3}{
                \node (s\i) [vertex, dotted] at (s.corner \i)
                {\TypeEsquaredlabel};
                \draw[->] (c) edge[bend right=8] (s\i) ;
                \draw[->] (s\i) edge[bend right=8, dotted] (c) ;
            }
            \draw[->] (s1) edge[out=-30,in=30, dotted, loop, looseness=5] (s1) ;
            \draw[->] (s3) edge[out=75,in=135, dotted, loop, looseness=5] (s3) ;

            \node (s6) [vertex, dotted] at (s.corner 6) {$A$};
            \draw[->] (c) edge[bend right=12] node[wt]{4} (s6) ;
            \draw[->] (s6) edge[bend right=12, dotted] (c) ;

            \foreach \i in {4,5,7,8}{
                \node (s\i) [vertex, dotted] at (s.corner \i) {$I$} ;
                \draw[->] (c) edge[bend right=12] node[wt]{2} (s\i) ;
                \draw[->] (s\i) edge[bend right=12, dotted] (c) ;
            }

            \draw[->] (c) edge[out=90,in=135,loop, looseness=10] (c) ;

            \draw[->] (s1) edge[bend right=30, dotted] (s3) ;
            \draw[->] (s3) edge[bend left=10, dotted] (s1) ;
        \end{tikzpicture}
    \end{tabular}
    \caption{%
        Neighbourhoods of the general \TypeEsquared and \TypeIII vertices.
    }
    \label{fig:TypeEsquared-TypeIII}
\end{figure}
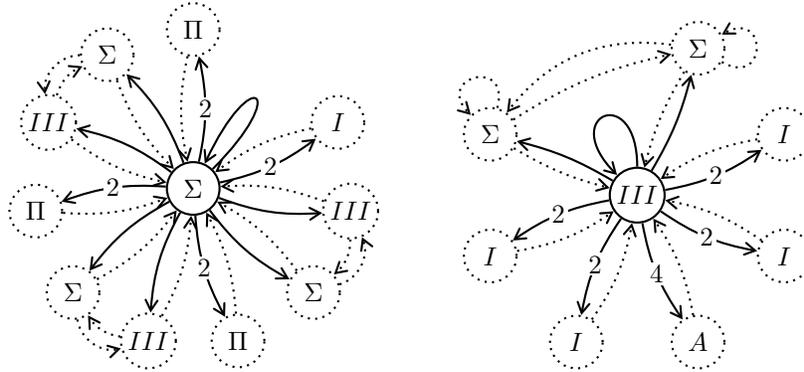

The general \TypeEsquared vertex is \(\classof{\EC\times\EC}\)
where \(\EC\) has no special automorphisms,
so
\(\RAut(\EC^2) = \subgrp{\sigma,\tau} \cong C_2^2\).
The orbits of the kernel subgroups
under \(\RAut(\EC^2)\) 
(with respect to an arbitrary labelling of \(\EC[2]\))
and the types of the corresponding neighbours
are described by Table~\ref{tab:TypeEsquared},
and
the neighbourhood of the generic \TypeEsquared vertex
is shown on the left of Figure~\ref{fig:TypeEsquared-TypeIII}.

\begin{table}[ht]
    \centering
    \begin{tabular}{c|c|c||c|c|c}
        Kernel orbit & Stab.  & Codomain 
        &
        Kernel orbit & Stab.  & Codomain 
        \\
        \hline
        \(\{K_{1,1}\}\) & \(\subgrp{\sigma,\tau}\) & \TypeEsquared
        &
        \(\{K_{\text{Id}}\}\) & \(\subgrp{\sigma,\tau}\) & (loop)
        \\
        \(\{K_{2,2}\}\) & \(\subgrp{\sigma,\tau}\) & \TypeEsquared
        &
        \(\{K_{(1,2)(3)}\}\) & \(\subgrp{\sigma,\tau}\) & \TypeIII
        \\
        \(\{K_{3,3}\}\) & \(\subgrp{\sigma,\tau}\) & \TypeEsquared
        &
        \(\{K_{(1,3)(2)}\}\) & \(\subgrp{\sigma,\tau}\) & \TypeIII
        \\
        \(\{K_{1,2}, K_{2,1}\}\) & \(\subgrp{\sigma}\) & \TypeExE
        &
        \(\{K_{(2,3)(1)}\}\) & \(\subgrp{\sigma,\tau}\) & \TypeIII
        \\
        \(\{K_{1,3}, K_{3,1}\}\) & \(\subgrp{\sigma}\) & \TypeExE
        &
        \(\{K_{(1,2,3)}, K_{(1,3,2)}\}\) & \(\subgrp{\sigma}\) & \TypeI
        \\
        \(\{K_{2,3}, K_{3,2}\}\) & \(\subgrp{\sigma}\) & \TypeExE
        \\
    \end{tabular}
    \caption{Edge data for the generic \TypeEsquared vertex.}
    \label{tab:TypeEsquared}
\end{table}

\subsection{\texorpdfstring{\TypeIII}{Type-III} vertices}
The generic \TypeIII vertex is 
\(\classof{\Jac{\XC_{III}}}\),
where
\[
    \XC_{III}: y^2 = (x^2-1)(x^2-u^2)(x^2-1/u^2)
\]
with \(u\) a free parameter;
note that \(\XC_{III}(u) = \XC_{I}(s,t)\)
with \((s,t) = (u,u^{-1})\).
We have
\(\RAut(\Jac{\XC}_{III}) = \subgrp{\sigma,\tau} \cong C_2^2\),
where \(\sigma\) is inherited from \TypeI
and \(\tau\) acts on \(x\)-coordinates via
\[
    \tau_*: x \longmapsto 1/x \,.
\]
Specializing the kernels and quadratic splittings
of~\S\ref{sec:TypeI}
at \((s,t) = (u,u^{-1})\),
we see that \(\RAut(\Jac{\XC_{III}})\)
acts on the kernel indices by
\begin{align*}
    \sigma_*
    & =
    (1)(2)(3)(4)(5)(6)(7)(8,9)(10,11)(12,13)(14,15)
    \,,
    \\
    \tau_*
    & =
    (1)(2)(3)(4,6)(5,7)(8)(9)(10,13)(11,12)(14)(15)
    \,.
\end{align*}
The kernel orbits and the edges leaving \(\classof{\Jac{\XC_{III}}}\)
are described in Table~\ref{tab:TypeIII}.

\begin{table}[ht]
    \centering
    \begin{tabular}{c|c|c||c|c|c}
        Orbit & Stab. & Codomain 
        &
        Orbit & Stab. & Codomain 
        \\
        \hline
        \(\{K_{1}\}\) & \(\subgrp{\sigma,\tau}\) & \TypeEsquared
        &
        \(\{K_{5},K_{7}\}\) & \(\subgrp{\sigma}\) & \TypeI
        \\
        \(\{K_{2}\}\) & \(\subgrp{\sigma,\tau}\) & (loop)
        &
        \(\{K_{8},K_{9}\}\) & \(\subgrp{\tau}\) & \TypeI
        \\
        \(\{K_{3}\}\) & \(\subgrp{\sigma,\tau}\) & \TypeEsquared
        &
        \(\{K_{i}: 10\le i\le 13\}\) & \(1\) & \TypeA
        \\
        \(\{K_{4},K_{6}\}\) & \(\subgrp{\sigma}\) & \TypeI
        &
        \(\{K_{14},K_{15}\}\) & \(\subgrp{\tau}\) & \TypeI
    \end{tabular}
    \caption{Edge data for the generic \TypeIII vertex.}
    \label{tab:TypeIII}
\end{table}

\begin{figure}[ht]
    \centering
    \begin{tikzpicture}[
            >={angle 60},
            thick,
            vertex/.style = {circle, draw, fill=white, inner sep=0.5mm, minimum size=6mm},
            wt/.style = {fill=white, anchor=center, pos=0.5, minimum size=1mm, inner sep=1pt}
        ]
        \node (s) [
            regular polygon,
            regular polygon sides=11,
            minimum size=60mm,
            rotate=-20,
            above
        ] at (0,0) {};
        \node (c) [vertex] at (s.center) {\TypeEsquaredlabel} ;
        \draw[->] (c) edge[out=54, in=82, loop, looseness=24] (c) ;
        \node (I) [vertex, dotted] at (s.corner 2) {$I$} ;
        \draw[->] (c) edge[bend right=6] node[wt]{2} (I) ;
        \draw[->] (I) edge[bend right=6,dotted] (c) ;

        \node (Sigma1) [vertex, dotted] at (s.corner 3) {\TypeEsquaredlabel} ;
        \node (III1) [vertex] at (s.corner 4) {$III$} ;
        \node (Pi1) [vertex, dotted] at (s.corner 5) {\TypeExElabel} ;
        
        \node (Sigma2) [vertex, dotted] at (s.corner 6) {\TypeEsquaredlabel} ;
        \node (III2) [vertex] at (s.corner 7) {$III$} ;
        \node (Pi2) [vertex, dotted] at (s.corner 8) {\TypeExElabel} ;
        
        \node (Sigma3) [vertex, dotted] at (s.corner 9) {\TypeEsquaredlabel} ;
        \node (III3) [vertex] at (s.corner 10) {$III$} ;
        \node (Pi3) [vertex, dotted] at (s.corner 11) {\TypeExElabel} ;

        \foreach \i in {1,2,3} {
            \draw[->] (c) edge[bend right=5] (Sigma\i) ;
            \draw[->] (Sigma\i) edge[bend right=5, dotted] (c) ;
            \draw[->] (c) edge[bend right=5] (III\i) ;
            \draw[->] (III\i) edge[bend right=5] (c) ;
            \draw[->] (Sigma\i) edge[bend right=12] (III\i) ;
            \draw[->] (III\i) edge[bend right=12,dotted] (Sigma\i) ;
            \draw[->] (Sigma\i) edge[out=100*\i-60, in=100*\i, dotted, loop, looseness=8] (Sigma\i) ;
            \draw[->] (III\i) edge[out=100*\i+190, in=100*\i+235, loop, looseness=8] (III\i) ;
            \node (IIIouter) [
                regular polygon,
                regular polygon sides=18,
                minimum size=48mm,
                shape border rotate=100*\i-40
            ] at (III\i) {};
            \foreach \j in {1,2,5,6}{
                \node (IIII\i-\j) [vertex, dotted] at (IIIouter.corner \j) {$I$} ;
                \draw[->] (III\i) edge[bend right=8] node[wt]{2} (IIII\i-\j) ;
                \draw[->] (IIII\i-\j) edge [bend right=8, dotted] (III\i) ;
            }
            \node (IIIA\i) [vertex, dotted] at (IIIouter.side 3) {$A$} ;
            \draw[->] (III\i) edge[bend right=8] node[wt]{4} (IIIA\i) ;
            \draw[->] (IIIA\i) edge [bend right=8, dotted] (III\i) ;

            \draw[->] (c) edge[bend right=6] node[wt]{2} (Pi\i) ;
            \draw[->] (Pi\i) edge[bend right=6, dotted] (c) ;
            \draw[->] (Pi\i) edge[bend right=6, dotted] (IIII\i-5) ;
            \draw[->] (IIII\i-5) edge[bend right=6, dotted] (Pi\i) ;
            \draw[->] (Pi\i) edge[bend right=6, dotted] (IIII\i-6) ;
            \draw[->] (IIII\i-6) edge[bend right=6, dotted] (Pi\i) ;

            \node (SigmaPi\i) [vertex, dotted] at (IIIouter.corner 17) {\TypeExElabel} ;
            \draw[->] (Sigma\i) edge[bend right=16, dotted] node[wt]{2} (SigmaPi\i) ;
            \draw[->] (SigmaPi\i) edge[bend right=16, dotted] (Sigma\i) ;
            \draw[->] (SigmaPi\i) edge[bend right=6, dotted] (IIII\i-1) ;
            \draw[->] (IIII\i-1) edge[bend right=6, dotted] (SigmaPi\i) ;
            \draw[->] (SigmaPi\i) edge[bend right=72, dotted] (IIII\i-2) ;
            \draw[->] (IIII\i-2) edge[bend right=-60, dotted] (SigmaPi\i) ;
        }
    \end{tikzpicture}
    \caption{%
        The neighbourhood of a generic \protect{\TypeEsquared}
        vertex and its \protect{\TypeIII} neighbours.
    }
    \label{fig:TypeEsquared-TypeIII-connect}
\end{figure}
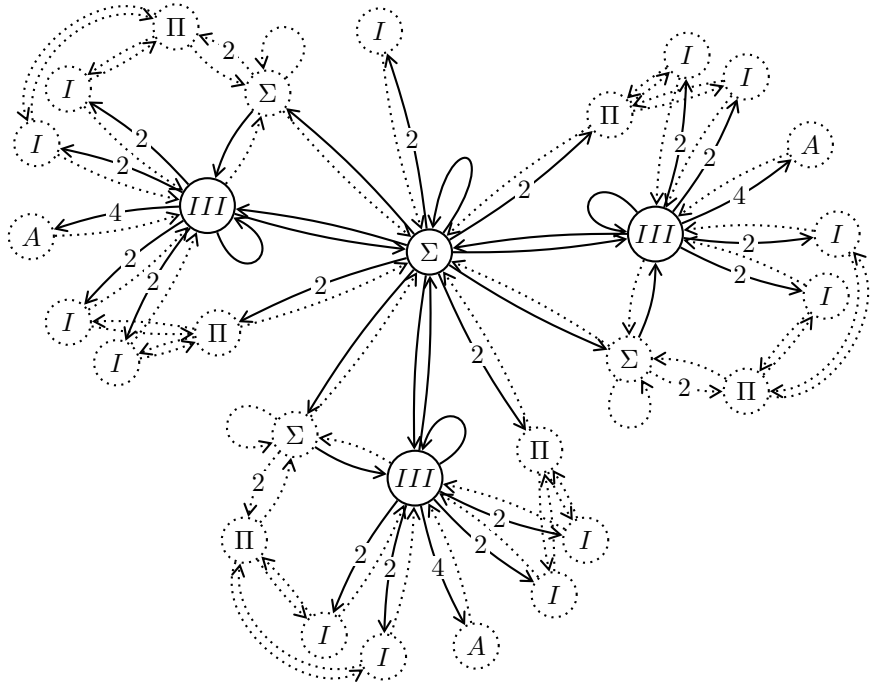

We observe that
\(\Jac{\XC_{III}}/K_2 \cong \Jac{\XC_{III}}\):
that is, \(\phi_{2}\) is a \((2,2)\)-endomorphism
of \(\Jac{\XC_{III}}\),
so \(\classof{\phi_{2}}\) is a weight-1 loop.
The kernels \(K_{1}\) and \(K_{3}\) 
are stabilised by \(\RAut(\Jac{\XC_{III}})\)
and \(\delta(K_{1}) = \delta(K_{3}) = 0\),
so \(\classof{\phi_{1}}\)
and \(\classof{\phi_{2}}\)
are weight-1 edges to \TypeEsquared
vertices
\(\classof{\EC^2}\)
and
\(\classof{(\EC')^2}\),
respectively,
where
\(\EC\) and \(\EC'\) are the elliptic curves
\begin{align*}
    \EC : y^2 & = (x - 1)(x - u^2)(x - 1/u^2)
    \,,
    \\
    \EC' : 
    y^2 
    & = 
    -2\left(x - 1\right)
        \Big(x^2 + 2\frac{u^4 - 6u^2 + 1}{(u^2+1)^2}x + 1\Big)
    \,.
\end{align*}
There is a \(2\)-isogeny
\(\varphi: \EC \to \EC'\),
as predicted in~\cite[\S4]{2001/Gaudry--Schost}
(in fact
\( \ker\varphi = \subgrp{(1,0)}\)
and
\( \ker\dualof{\varphi} = \subgrp{(1,0)}\)),
so there are edges \(\classof{\varphi\times\varphi}\)
and \(\classof{\dualof{\varphi}\times\dualof{\varphi}}\)
between \(\classof{\EC^2}\)
and \(\classof{(\EC')^2}\).
The neighbourhood of the general \TypeIII vertex
is shown on the right of Figure~\ref{fig:TypeEsquared-TypeIII}.
Combining with the \TypeEsquared neighbourhood
and extending to include shared adjacent vertices
yields Figure~\ref{fig:TypeEsquared-TypeIII-connect}.

\subsection{Elliptic 3-isogenies: \texorpdfstring{\TypeIV}{Type-IV} vertices}

\begin{figure}[ht]
    \centering
    \begin{tikzpicture}[
            >={angle 60},
            thick,
            vertex/.style = {circle, draw, fill=white, inner sep=0.5mm, minimum size=7mm},
            wt/.style = {fill=white, anchor=center, pos=0.5, minimum size=1mm, inner sep=1pt}
        ]
        \node (s) [
            regular polygon,
            regular polygon sides=7,
            minimum size=42mm,
            above
        ] at (0,0) {};
        \node (c) [vertex] at (s.center) {$IV$};
        \node (s1) [vertex, dotted] at (s.corner 1)
        {\TypePhilabel} ;
        \draw[->] (c) edge[bend right=12] node[wt]{3} (s1) ;
        \draw[->] (s1) edge[bend right=12, dotted] (c) ;

        \foreach \i in {2,4,6}{
            \node (s\i) [vertex, dotted] at (s.corner \i) {$I$} ;
            \draw[->] (c) edge[bend right=12] node[wt]{3} (s\i) ;
            \draw[->] (s\i) edge[bend right=12, dotted] (c) ;
        }
        \foreach \i in {3,5,7}{
            \node (s\i) [vertex, dotted] at (s.corner \i) {$IV$} ;
            \draw[->] (c) edge[bend right=12] (s\i) ;
            \draw[->] (s\i) edge[bend right=12, dotted] (c) ;
        }
    \end{tikzpicture}
    \caption{The neighbourhood of the general \TypeIV vertex.}
    \label{fig:TypeIV}
\end{figure}
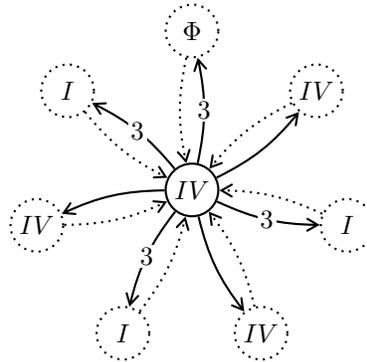

The generic \TypeIV vertex is represented by 
\(\Jac{\XC_{IV}(v)}\),
where \(\XC_{IV}(v) := \XC_{I}(s_{IV}(v),t_{IV}(v))\)
with
\[
    s_{IV}(v) 
    :=
    \frac{
        (v+1)(v-\zeta_3)
    }{
        (v-1)(v+\zeta_3)
    }
    \qquad\text{and}\qquad
    t_{IV}(v) 
    := 
    \frac{(v+1)(v-\zeta_3^2)}{(v-1)(v+\zeta_3^2)}
\]
where \(\zeta_3\) is a primitive third root of unity
and \(v\) is a free parameter.
We have
\(\RAut(\XC_{IV}(v)) = \subgrp{\sigma,\rho}\cong S_3\),
where \(\sigma\) is inherited from \TypeI
and \(\rho\) is the order-3 automorphism acting 
on \(x\)-coordinates via
\[
    \rho_*: 
    x 
    \longmapsto
    \big( (2\zeta_3+1)(v^2-1)x + 3(v+1)^2 \big)
    /
    \big( 3(v-1)^2x + (2\zeta_3+1)(v^2-1) \big)
    \,.
\]
Specializing the kernels and quadratic splittings
from~\S\ref{sec:TypeI},
we see that the action of \(\RAut(\Jac{\XC_{IV}})\) 
on (the indices of) the \(K_i\)
is given by
\begin{align*}
    \rho_*
    & = 
    (1,9,8)
    (2,15,14)
    (3)
    (4,12,13)
    (5)
    (6,10,11)
    (7)
    \,,
    \\
    \sigma_*
    & =
    (1)(2)(3)(4)(5)(6)(7)(8,9)(10,11)(12,13)(14,15)
    \,.
\end{align*}
The kernel orbits
and the edges leaving \(\classof{\Jac{\XC}_{IV}}\)
described in Table~\ref{tab:TypeIV},
and illustrated in Figure~\ref{fig:TypeIV}.
We find that
\(\classof{\AV_{1}} = \classof{\AV_{8}} = \classof{\AV_{9}}
= \classof{\EC\times\EC'}\),
where
\begin{align*}
    \EC: y^2 & = (x-1)(x-s_{IV}(v)^2)(x - t_{IV}(v)^2)
    \shortintertext{and}
    \EC': y^2 & = (x-1)(x-1/s_{IV}(v)^2)(x - 1/t_{IV}(v)^2)
    \,.
\end{align*}
There is a \(3\)-isogeny
\(\Phi: \EC \to \EC'\),
as predicted in~\cite[\S3]{2001/Gaudry--Schost}:
the kernel of \(\Phi\)
is cut out by \(x - (v+1)^2/(v-1)^2\),
and the kernel of \(\dualof{\Phi}\)
is cut out by \(x - (v-1)^2/(v+1)^2\).

Elliptic products with a 3-isogeny between the factors
therefore play a special role in the Richelot isogeny graph;
we will represent these special \TypeExE vertices
using the symbol \(\Phi\).
We remark that the presence of the 3-isogeny
severely constrains the possible specializations
of a \(\Phi\)-vertex.

\begin{table}[ht]
    \centering
    \begin{tabular}{c|c|c||c|c|c}
        Kernel orbit
        & Stabilizer
        & Codomain
        & Kernel orbit
        & Stabilizer
        & Codomain
        \\
        \hline
        \(\{K_{1}, K_{8}, K_{9}\}\) 
        & \(\subgrp{\sigma}\)
        & \TypeExE (\TypePhilabel)
        &
        \(\{K_{3}\}\)
        & \(\subgrp{\sigma,\rho}\)
        & \TypeIV
        \\
        \(\{K_{2}, K_{14}, K_{15}\}\) 
        & \(\subgrp{\sigma}\)
        & \TypeI
        &
        \(\{K_{5}\}\)
        & \(\subgrp{\sigma,\rho}\)
        & \TypeIV
        \\
        \(\{K_{4}, K_{12}, K_{13}\}\) 
        & \(\subgrp{\sigma}\)
        & \TypeI
        &
        \(\{K_{7}\}\)
        & \(\subgrp{\sigma,\rho}\)
        & \TypeIV
        \\
        \(\{K_{6}, K_{10}, K_{11}\}\) 
        & \(\subgrp{\sigma}\)
        & \TypeI
        \\
    \end{tabular}
    \caption{Edge data for the generic \TypeIV vertex.}
    \label{tab:TypeIV}
\end{table}

\begin{figure}[ht]
    \centering
    \begin{tikzpicture}[
            >={angle 60},
            thick,
            vertex/.style = {circle, draw, fill=white, inner sep=0.5mm, minimum size=6mm},
            wt/.style = {fill=white, anchor=center, pos=0.5, minimum size=1mm, inner sep=1pt}
        ]
        \node (domain) [vertex, ultra thick] at (0,0) {$IV$} ;
        \node (codomain) [vertex, ultra thick] at (7,0) {\TypePhilabel} ;
        \draw[->, ultra thick] (codomain) edge[bend right=6] (domain) ;
        \draw[->, ultra thick] (domain) edge[bend right=6] node[wt]{3} (codomain) ;

        \node (Xtop) [
            regular polygon,
            regular polygon sides=4,
            minimum size=22mm,
        ] at (3.5, 3.8) {};
        \node (Xtopcod1) [vertex, dotted] at (Xtop.corner 1) {\TypePhilabel} ;
        \node (Xtopdom1) [vertex, dotted] at (Xtop.corner 2) {$IV$} ;
        \node (Xtopdom2) [vertex, dotted] at (Xtop.corner 3) {$I$} ;
        \node (Xtopcod2) [vertex, dotted] at (Xtop.corner 4) {$I$} ;
        \foreach \i in {1,2} {
            \draw[->] (Xtopdom\i) edge[bend right=40, dotted] (domain) ;
            \draw[->] (codomain) edge[bend right=48] (Xtopcod\i) ;
            \draw[->] (Xtopcod\i) edge[bend right=-40, dotted] (codomain) ;
            \foreach \j in {1,2} {
                \draw[->] (Xtopcod\j) edge[bend right=14, dotted] (Xtopdom\i) ;
            }
        }
        \draw[->] (Xtopdom1) edge[bend right=14, dotted] node[wt]{3} (Xtopcod1) ;
        \draw[->] (Xtopdom1) edge[bend right=14, dotted] node[wt]{3} (Xtopcod2) ;
        \draw[->] (Xtopdom2) edge[bend right=14, dotted] (Xtopcod1) ;
        \draw[->] (Xtopdom2) edge[bend right=14, dotted] (Xtopcod2) ;
        \draw[->] (domain) edge[bend right=-48] (Xtopdom1) ;
        \draw[->] (domain) edge[bend right=-52] node[wt]{3} (Xtopdom2) ;

        \node (Xmid) [
            regular polygon,
            regular polygon sides=4,
            minimum size=22mm,
        ] at (3.5, 1.5) {};
        \node (Xmidcod1) [vertex, dotted] at (Xmid.corner 1) {\TypePhilabel} ;
        \node (Xmiddom1) [vertex, dotted] at (Xmid.corner 2) {$IV$} ;
        \node (Xmiddom2) [vertex, dotted] at (Xmid.corner 3) {$I$} ;
        \node (Xmidcod2) [vertex, dotted] at (Xmid.corner 4) {$I$} ;
        \foreach \i in {1,2} {
            \draw[->] (Xmiddom\i) edge[bend right=16, dotted] (domain) ;
            \draw[->] (codomain) edge[bend right=24] (Xmidcod\i) ;
            \draw[->] (Xmidcod\i) edge[bend right=-16, dotted] (codomain) ;
            \foreach \j in {1,2} {
                \draw[->] (Xmidcod\j) edge[bend right=14, dotted] (Xmiddom\i) ;
            }
        }
        \draw[->] (Xmiddom1) edge[bend right=14, dotted] node[wt]{3} (Xmidcod1) ;
        \draw[->] (Xmiddom1) edge[bend right=14, dotted] node[wt]{3} (Xmidcod2) ;
        \draw[->] (Xmiddom2) edge[bend right=14, dotted] (Xmidcod1) ;
        \draw[->] (Xmiddom2) edge[bend right=14, dotted] (Xmidcod2) ;
        \draw[->] (domain) edge[bend right=-24] (Xmiddom1) ;
        \draw[->] (domain) edge[bend right=-28] node[wt]{3} (Xmiddom2) ;

        \node (Xbot) [
            regular polygon,
            regular polygon sides=4,
            minimum size=22mm,
        ] at (3.5, -1.5) {};
        \node (Xbotcod1) [vertex, dotted] at (Xbot.corner 1) {\TypePhilabel} ;
        \node (Xbotdom1) [vertex, dotted] at (Xbot.corner 2) {$IV$} ;
        \node (Xbotdom2) [vertex, dotted] at (Xbot.corner 3) {$I$} ;
        \node (Xbotcod2) [vertex, dotted] at (Xbot.corner 4) {$I$} ;
        \foreach \i in {1,2} {
            \draw[->] (Xbotdom\i) edge[bend right=-16, dotted] (domain) ;
            \draw[->] (codomain) edge[bend right=-24] (Xbotcod\i) ;
            \draw[->] (Xbotcod\i) edge[bend right=16, dotted] (codomain) ;
            \foreach \j in {1,2} {
                \draw[->] (Xbotcod\j) edge[bend right=14, dotted] (Xbotdom\i) ;
            }
        }
        \draw[->] (Xbotdom1) edge[bend right=14, dotted] node[wt]{3} (Xbotcod1) ;
        \draw[->] (Xbotdom1) edge[bend right=14, dotted] node[wt]{3} (Xbotcod2) ;
        \draw[->] (Xbotdom2) edge[bend right=14, dotted] (Xbotcod1) ;
        \draw[->] (Xbotdom2) edge[bend right=14, dotted] (Xbotcod2) ;
        \draw[->] (domain) edge[bend right=24] (Xbotdom1) ;
        \draw[->] (domain) edge[bend right=28] node[wt]{3} (Xbotdom2) ;

        \node (codomainout) [
            regular polygon,
            regular polygon sides=20,
            minimum size=42mm,
        ] at (codomain) {};
        \foreach \i in {13,20} {
            \node (codA\i) [vertex, dotted] at (codomainout.corner \i) {$I$} ;
            \draw[->] (codomain) edge[bend right=4] (codA\i) ;
            \draw[->] (codA\i) edge[bend right=4, dotted] (codomain) ;
        }
        \foreach \i in {14,15,16,17,18,19} {
            \node (codA\i) [vertex, dotted] at (codomainout.corner \i) {\TypeExElabel} ;
            \draw[->] (codomain) edge[bend right=4] (codA\i) ;
            \draw[->] (codA\i) edge[bend right=4, dotted] (codomain) ;
        }

    \end{tikzpicture}
    \caption{%
        The neighbourhood of a
        \TypeIV vertex and its \TypeExE neighbour.  
    }
    \label{fig:TypeIV-TypeExE-connect}
\end{figure}
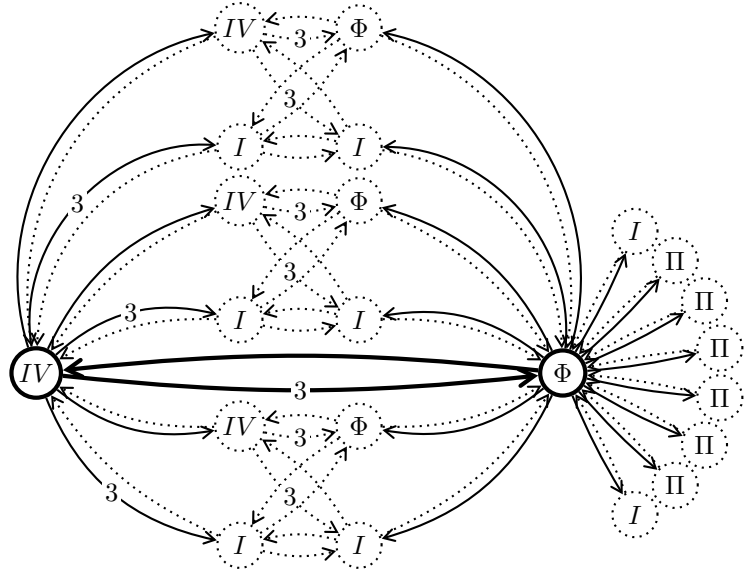

Figure~\ref{fig:TypeIV-TypeExE-connect}
shows the neighbourhood of the edges between a general \TypeIV vertex
and its \TypePhilabel-neighbour;
it should be compared with Figures~\ref{fig:general-edge}
and~\ref{fig:TypeI-TypeExE-connect}.
\TypeIV vertices correspond to \TypePhilabel-vertices,
and edges between \TypeIV vertices
correspond to edges between \TypePhilabel-vertices.

\subsection{The \texorpdfstring{\TypeExEzero}{Type-Pi-0} family}
\label{sec:TypeExEzero}

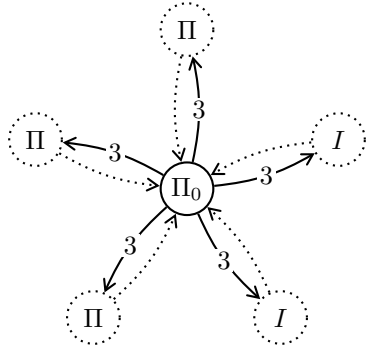
\begin{figure}[ht]
    \begin{tikzpicture}[
            >={angle 60},
            thick,
            vertex/.style = {circle, draw, fill=white, inner sep=0.5mm, minimum size=7mm},
            wt/.style = {fill=white, anchor=center, pos=0.5, minimum size=1mm, inner sep=1pt}
        ]
        \node (s) [
            regular polygon,
            regular polygon sides=5,
            minimum size=42mm,
            above] at (0,0) {};
        \node (c) [vertex] at (s.center) {\TypeExEzerolabel};
        \foreach \i in {1,2,3}{
            \node (s\i) [vertex, dotted] at (s.corner \i)
            {\TypeExElabel} ;
            \draw[->] (c) edge[bend right=12] node[wt]{3} (s\i) ;
            \draw[->] (s\i) edge[bend right=12, dotted] (c) ;
        }
        \foreach \i in {4,5}{
            \node (s\i) [vertex, dotted] at (s.corner \i) {$I$} ;
            \draw[->] (c) edge[bend right=12] node[wt]{3} (s\i) ;
            \draw[->] (s\i) edge[bend right=12, dotted] (c) ;
        }
    \end{tikzpicture}
    \caption{The neighbourhood of a generic \TypeExEzero vertex.}
    \label{fig:TypeExEzero}
\end{figure}

The \TypeExEzero vertices
are \(\classof{\EC\times\EC_{0}}\) 
for elliptic curves \(\EC\not\cong\EC_{0}\).
We have 
\(
    \RAut(\EC\times\EC_{0}) 
    = 
    \subgrp{\sigma,[1]\times\zeta} 
    \cong 
    C_6
\).
The automorphism \(\zeta\) of \(\EC_0\)
cycles the points of order 2 on \(\EC_0\),
so \([1]\times\zeta\) fixes no \((2,2)\)-isogeny kernels.
Instead, the kernel subgroups of \(\EC\times\EC_0[2]\)
form orbits of three, 
and so we see the five neighbours with weight-3 edges
in Figure~\ref{fig:TypeExEzero}
(which should be compared with Figure~\ref{fig:TypeExE-TypeI}).

\subsection{The \texorpdfstring{\TypeExEtwelvecubed}{Type-Pi-1728} family}
\label{sec:TypeExEtwelvecubed}

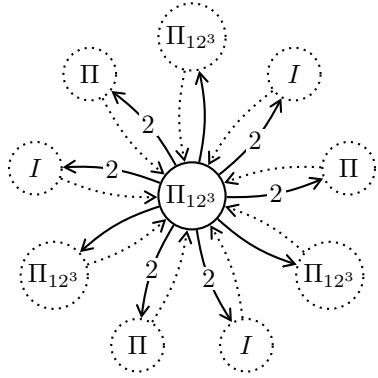
\begin{figure}[ht]
    \begin{tikzpicture}[
            >={angle 60},
            thick,
            vertex/.style = {circle, draw, fill=white, inner sep=0.5mm, minimum size=7mm},
            wt/.style = {fill=white, anchor=center, pos=0.5, minimum size=1mm, inner sep=1pt}
        ]
        \node (s) [
            regular polygon,
            regular polygon sides=9,
            minimum size=42mm,
            above
        ] at (0,0) {};
        \node (c) [vertex] at (s.center) {\TypeExEtwelvecubedlabel};
        
        \foreach \i in {1,4,7}{
            \node (s\i) [vertex, dotted] at (s.corner \i) {\TypeExEtwelvecubedlabel} ;
            \draw[->] (c) edge[bend right=12] (s\i) ;
            \draw[->] (s\i) edge[bend right=12, dotted] (c) ;
        }
        \foreach \i in {2,5,8}{
            \node (s\i) [vertex, dotted] at (s.corner \i)
            {\TypeExElabel} ;
            \draw[->] (c) edge[bend right=12] node[wt]{2} (s\i) ;
            \draw[->] (s\i) edge[bend right=12, dotted] (c) ;
        }
        \foreach \i in {3,6,9}{
            \node (s\i) [vertex, dotted] at (s.corner \i) {$I$} ;
            \draw[->] (c) edge[bend right=12] node[wt]{2} (s\i) ;
            \draw[->] (s\i) edge[bend right=12, dotted] (c) ;
        }
    \end{tikzpicture}
    \caption{%
        The neighbourhoods of the generic \TypeExEtwelvecubed vertex.
    }
    \label{fig:TypeExEtwelvecubed}
\end{figure}

The \TypeExEtwelvecubed vertices
are \(\classof{\EC\times\EC_{12^3}}\)
for elliptic curves \(\EC\not\cong\EC_{12^3}\).
The curve \(\EC_{12^3}\) has an order-4 automorphism \(\iota\)
which fixes one point \(P_3\) of order 2,
and exchanges \(P_1\) and \(P_2\).
We therefore have 
an order-4 element \(\alpha = [1]\times\iota\) 
generating \(\RAut(\EC\times\EC_{12^3}) \cong C_4\),
and \(\alpha^2 = [1]\times[-1] = \sigma\) (which fixes all the kernels).
Hence, with respect to \(\subgrp{\alpha}\),
the isometries form three orbits of size two,
as do the six product kernels not involving \(P_3\);
on the other hand, the kernels \(\subgrp{P}\times\subgrp{P_3}\)
are fixed by \(\alpha\),
and since \(\EC_{12^3}/\subgrp{P_1}\cong\EC_{12^3}\)
we get three weight-1 edges
to \TypeExEtwelvecubed vertices.
The situation is illustrated in Figure~\ref{fig:TypeExEtwelvecubed}
(which should be compared with Figure~\ref{fig:TypeExE-TypeI}).

\subsection{The \texorpdfstring{\TypeEzeroxEtwelvecubed}{Type-Pi-0-1728} vertex}

The unique \TypeEzeroxEtwelvecubed vertex is
\(\classof{\EC_0\times\EC_{12^3}}\).
Its reduced automorphism group is
\(
    \RAut(\EC_0\times\EC_{12^3})
    = 
    \subgrp{\zeta\times[1], [1]\times\iota}
    \cong C_{12}
\).
The kernel orbits and edges 
can be derived using a combination of the analyses
in~\S\ref{sec:TypeExEzero}
and~\S\ref{sec:TypeExEtwelvecubed};
the results
are described in Table~\ref{tab:TypeEzeroxEtwelvecubed}.
The neighbourhood of the \TypeEzeroxEtwelvecubed vertex,
illustrated in Figure~\ref{fig:TypeEzeroxEtwelvecubed},
is a combination of the \TypeExEzero
and \TypeExEtwelvecubed neighbourhoods
of Figures~\ref{fig:TypeExEzero} and~\ref{fig:TypeExEtwelvecubed}.

\begin{table}[ht]
    \centering
    \begin{tabular}{c|c|ccc}
        \multirow{2}{*}{Kernel orbit}
        & Stabilizer 
        & \multicolumn{3}{c}{Codomain type}
        \\
        & (conjugate)
        & General \(p\)
        & \(p = 7\)
        & \(p = 11\)
        \\
        \hline
        \(\{K_{1,3},K_{2,3},K_{3,3}\}\) 
            & \(\subgrp{[1]\times\iota}\)
            & \TypeExEtwelvecubed
            & \TypeExEtwelvecubed
            & \TypeEtwelvecubedsquared
        \\
        \(\{K_{i,j} : 1 \le i \le 3, 1 \le j \le 2\}\)
            & \(\subgrp{\sigma}\)
            & \TypeExE
            & \TypeExEtwelvecubed
            & (loops)
        \\
        \(\{K_{\pi}: \pi \in S_3\}\) 
            & \(\subgrp{\sigma}\)
            & \TypeI
            & \TypeI
            & \TypeIV
        \\
    \end{tabular}
    \caption{Edge data for the unique \TypeEzeroxEtwelvecubed vertex.}
    \label{tab:TypeEzeroxEtwelvecubed}
\end{table}

\begin{figure}[ht]
    \centering
    \begin{tikzpicture}[
            >={angle 60},
            thick,
            vertex/.style = {circle, draw, fill=white, inner sep=0.5mm, minimum size=7mm},
            wt/.style = {fill=white, anchor=center, pos=0.5, minimum size=1mm, inner sep=1pt}
        ]
        \node (c) [vertex] at (0,0) {\TypeEzeroxEtwelvecubedlabel};
        \node (s1) [vertex, dotted] at (0,2) {\TypeExEtwelvecubedlabel} ;
        \draw[->] (c) edge[bend right=12] node[wt]{3} (s1) ;
        \draw[->] (s1) edge[bend right=12, dotted] (c) ;

        \node (s2) [vertex, dotted] at (-2,0)
        {\TypeExElabel} ;
        \draw[->] (c) edge[bend right=12] node[wt]{6} (s2) ;
        \draw[->] (s2) edge[bend right=12, dotted] (c) ;

        \node (s3) [vertex, dotted] at (2,0) {$I$} ;
        \draw[->] (c) edge[bend right=12] node[wt]{6} (s3) ;
        \draw[->] (s3) edge[bend right=12, dotted] (c) ;
    \end{tikzpicture}
    \caption{%
        The neighbourhood of the \TypeEzeroxEtwelvecubed vertex.
    }
    \label{fig:TypeEzeroxEtwelvecubed}
\end{figure}

\subsection{The \texorpdfstring{\TypeEtwelvecubedsquared}{Type-Sigma-1728} vertex}

The unique \TypeEtwelvecubedsquared vertex is
\(\classof{\EC_{12^3}^2}\).
We have
\(
    \RAut(\EC_{12^3}^2)
    = 
    \subgrp{\sigma,\tau,[1]\times\iota}
    \cong C_2^2\rtimes C_4
\).
The kernel orbits and edges are described in
Table~\ref{tab:TypeEtwelvecubedsquared}.
Figure~\ref{fig:TypeEtwelvecubedsquared},
illustrating the neighbourhood of \(\classof{\EC_{12^3}^2}\),
should be compared with Figure~\ref{fig:TypeEsquared-TypeIII}.

\begin{figure}[ht]
    \centering
    \begin{tikzpicture}[
            >={angle 60},
            thick,
            vertex/.style = {circle, draw, fill=white, inner sep=0.5mm, minimum size=7mm},
            wt/.style = {fill=white, anchor=center, pos=0.5, minimum size=1mm, inner sep=1pt}
        ]
        \node (c) [vertex] at (0,0) {\TypeEtwelvecubedsquaredlabel};
        \draw[->] (c) edge[out=-10,in=40,loop, looseness=6] (c) ;
        \draw[->] (c) edge[out=250,in=300,loop, looseness=6] node[wt]{2} (c) ;
        \node (s2) [vertex, dotted] at (-4,0) {\TypeEsquaredlabel} ;
        \draw[->] (c) edge[bend right=8] node[wt]{4} (s2) ;
        \draw[->] (s2) edge[bend right=8, dotted] (c) ;
        \draw[->] (s2) edge[loop, out=150, in=210, dotted, looseness=6] (s2) ;
        \node (s3) [vertex, dotted] at (-2,-2) {$III$} ;
        \draw[->] (c) edge[bend right=8] node[wt]{4} (s3) ;
        \draw[->] (s3) edge[bend right=8, dotted] (c) ;
        \draw[->] (s3) edge[loop, out=160, in=220, dotted, looseness=6] (s3) ;

        \draw[->] (s2) edge[bend right=8, dotted] (s3) ;
        \draw[->] (s3) edge[bend right=8, dotted] (s2) ;

        \node (s4) [vertex, dotted] at (3,-2) {\TypeExEtwelvecubedlabel} ;
        \draw[->] (c) edge[bend right=8] node[wt]{4} (s4) ;
        \draw[->] (s4) edge[bend right=8, dotted] (c) ;

    \end{tikzpicture}
    \caption{%
        The neighbourhood of the \TypeEtwelvecubedsquared vertex.
        The dotted neighbour types change for \(p = 7\) and \(11\)
        (see Table~\ref{tab:TypeEtwelvecubedsquared}).
    }
    \label{fig:TypeEtwelvecubedsquared}
\end{figure}

\begin{table}[ht]
    \centering
    \begin{tabular}{c|c|c@{\;}c@{\;}c}
        \multirow{2}{*}{Kernel orbit}
        & Stabilizer
        & \multicolumn{3}{c}{Codomain type}
        \\
        & (conjugate)
        & General \(p\)
        & \(p = 7\)
        & \(p = 11\)
        \\
        \hline
        \(\{K_{1,1},K_{1,2},K_{2,1},K_{2,2}\}\) 
            & \(\subgrp{\sigma,\tau}\) 
            & \TypeEsquared
            & (loops)
            & \TypeEzerosquared
        \\
        \(\{K_{1,3},K_{2,3},K_{3,1},K_{3,2}\}\) 
            & \(\subgrp{[1]\times\iota}\) 
            & \TypeExEtwelvecubed
            & (loops)
            & \TypeEzeroxEtwelvecubed
        \\
        \(\{K_{3,3}\}\) 
            & \(\RAut(\EC_{12^3})\) 
            & (loop)
            & (loop)
            & (loop)
        \\
        \(\{K_{\text{Id}},K_{(1,2)(3)}\}\) 
            & \(\subgrp{\tau,\iota\times\iota}\) 
            & (loops)
            & (loops)
            & (loops)
        \\
        \(\left\{\!\!\!\!\begin{array}{c}K_{(1,2,3)}, K_{(1,3,2)},\\
                                 K_{(1,3)(2)},
        K_{(1)(2,3)}\end{array}\!\!\!\!\right\}\)
            & \(\subgrp{\sigma,\tau}\) 
            & \TypeIII
            & \TypeIII
            & \TypeV
        \\
    \end{tabular}
    \caption{Edge data for the unique \TypeEtwelvecubedsquared vertex.}
    \label{tab:TypeEtwelvecubedsquared}
\end{table}

\subsection{The \texorpdfstring{\TypeV}{Type-V} and \texorpdfstring{\TypeEzerosquared}{Type-Sigma-0} vertices}

The \TypeV and \TypeEzerosquared vertices are always neighbours,
so we treat them simultaneously.

The unique \TypeEzerosquared vertex is \(\classof{\EC_{0}^2}\),
and \(\RAut(\EC_{0}^2) =
\subgrp{\tau,[1]\times\zeta,\zeta\times[1]}/\subgrp{-1} \cong C_6\times S_3\).
The kernel orbits and edges are described in
Table~\ref{tab:TypeEzerosquared}.

\begin{table}[ht]
    \centering
    \begin{tabular}{c|c|cc}
        \multirow{2}{*}{Kernel orbit}
        & Stabilizer
        & \multicolumn{2}{c}{Codomain type}
        \\
        & (conjugate)
        & General \(p\)
        & \(p = 11\)
        \\
        \hline
        \(\{K_{i,j} : 1 \le i, j \le 3\}\) 
            & \(\subgrp{\sigma,\tau}\) 
            & \TypeEsquared
            & \TypeEtwelvecubedsquared
        \\
        \(\{K_{\text{Id}}, K_{(1,2,3)}, K_{(1,3,2)}\}\)
            & \(\subgrp{\tau,\zeta\times(-\zeta)}\)
            & (loop)
            & (loop)
        \\
        \(\{K_{(1,2)(3)}, K_{(1,3)(2)}, K_{(2,3)(1)}\}\)
            & \(\subgrp{\tau, \zeta\times\zeta^2}\)
            & \TypeV
            & \TypeV
        \\
    \end{tabular}
    \caption{Edge data for the unique \TypeEzerosquared vertex.}
    \label{tab:TypeEzerosquared}
\end{table}

The unique \TypeV vertex is \(\classof{\Jac{\XC_V}}\),
where \(\XC_V: y^2 = x^6 + 1\);
note that \(\XC_{V} = \XC_{III}(\zeta_6) = \XC_{I}(\zeta_6,1/\zeta_6)\),
where \(\zeta_6\) is a primitive sixth root of unity.
We have \(\RAut(\XC_{V}) = \subgrp{\sigma,\tau,\zeta}\),
where \(\sigma\) and \(\tau\) are inherited from \(\XC_{III}\),
and \(\zeta\) is a new automorphism of order \(6\)
such that \(\zeta^3 = \sigma\).
Specializing the kernels and quadratic splittings
from~\S\ref{sec:TypeI},
these automorphisms act on (the indices of) the \(K_i\)
via 
\begin{align*}
    \tau_*
    & =
    (1)(2)(3)(4,6)(5,7)(8,9)(10,13)(11,12)(14,15)
    \,,
    \\
    \zeta_* 
    & = 
    (1)(2,4,6)(3,5,7)(8,9)(10,14,12,11,15,13)
    \,,
    \\
    \sigma_* = \zeta_*^3
    & =
    (1)(2)(3)(4)(5)(6)(7)(8,9)(10,11)(12,13)(14,15)
    \,.
\end{align*}
The kernel orbits and edges are described in Table~\ref{tab:TypeV}.

Figure~\ref{fig:TypeV} illustrates the shared neighbourhood 
of the \TypeV and \TypeEzerosquared vertices
for general \(p\);
it should be compared with Figure~\ref{fig:TypeEsquared-TypeIII}.
The \TypeI neighbour of the \TypeV vertex always has four \((2,2)\)-endomorphisms,
and they are included here for completeness,
as well as
the \TypeI and \TypeExE neighbours of the \TypeIV vertex,
since these are also connected to the \TypeEsquared and \TypeI
neighbours.
Dashed neighbour types may change for \(p = 11\), \(17\), \(29\), and \(41\)
(see Table~\ref{tab:TypeV}).

\begin{figure}[ht]
    \centering
    \begin{tikzpicture}[
            >={angle 60},
            thick,
            vertex/.style = {circle, draw, fill=white, inner sep=0.5mm, minimum size=7mm},
            wt/.style = {fill=white, anchor=center, pos=0.5, minimum size=1mm, inner sep=1pt}
        ]
        \node (cV) [vertex] at (0,0) {$V$};
        \node (ezerosquared) [vertex] at (-4,0) {\TypeEzerosquaredlabel};
        \node (cI) [vertex, dotted] at (4,0) {$I$};
        \node (esquared) [vertex, dotted] at (-2,-2) {\TypeEsquaredlabel};
        \node (cIV) [vertex, dotted] at (1,-1.3) {$IV$};
        \node (ee) [vertex, dotted] at (3,-1.3) {\TypePhilabel} ;
        \node (cInew) [vertex, dotted] at (4,-2) {$I$} ;

        \draw[->] (cV) edge[bend right=8] (ezerosquared) ;
        \draw[->] (ezerosquared) edge[bend right=8] node[wt]{3} (cV) ;
        \draw[->] (ezerosquared) edge[out=140,in=220,loop, looseness=5] node[wt]{3} (ezerosquared) ;
        \draw[->] (cV) edge[out=50,in=130,loop, looseness=5] node[wt]{3} (cV) ;

        \draw[->] (cV) edge[bend right=8] node[wt]{3} (esquared) ;
        \draw[->] (esquared) edge[bend right=8, dotted] (cV) ;
        \draw[->] (esquared) edge[out=-190,in=-110,loop, looseness=5] (esquared) ;

        \draw[->] (cV) edge[bend right=12] node[wt]{2} (cIV) ;
        \draw[->] (cIV) edge[bend right=12, dotted] (cV) ;

        \draw[->] (cV) edge[bend right=8] node[wt]{6} (cI) ;
        \draw[->] (cI) edge[bend right=8, dotted] (cV) ;
        \draw[->] (cI) edge[out=-40,in=40,loop, looseness=5] node[wt]{4} (cI) ;

        \draw[->] (ezerosquared) edge[bend right=8] node[wt]{9} (esquared) ;
        \draw[->] (esquared) edge[bend right=8, dotted] (ezerosquared) ;

        \draw[->] (esquared) edge[bend right=16, dotted] node[wt]{2} (ee) ;
        \draw[->] (ee) edge[bend right=-8, dotted] (esquared) ;
        \draw[->] (cIV) edge[bend right=-8, dotted] node[wt]{3} (ee) ;
        \draw[->] (ee) edge[bend right=20, dotted] (cIV) ;
        \draw[->] (cI) edge[bend right=12, dotted] (ee) ;
        \draw[->] (ee) edge[bend right=12, dotted] (cI) ;

        \draw[->] (esquared) edge[bend right=20, dotted] node[wt]{3} (cInew) ;
        \draw[->] (cInew) edge[bend right=-12, dotted] (esquared) ;
        \draw[->] (cIV) edge[bend right=20, dotted] node[wt]{3} (cInew) ;
        \draw[->] (cInew) edge[bend right=-8, dotted] (cIV) ;
        \draw[->] (cI) edge[bend right=8, dotted] (cInew) ;
        \draw[->] (cInew) edge[bend right=8, dotted] (cI) ;
 
    \end{tikzpicture}
    \caption{%
        The neighbourhood of the \protect{\TypeV} 
        and \protect{\TypeEzerosquared} vertices.
    }
    \label{fig:TypeV}
\end{figure}

\begin{table}[ht]
    \centering
    \begin{tabular}{c|c|c@{\;}c@{\;}c@{\;}c@{\;}c}
        \multirow{2}{*}{Kernel orbit}
        & Stabilizer
        & \multicolumn{5}{c}{Codomain type} 
        \\
        & (conj.)
        & General \(p\)
        & \(p = 11\)
        & \(p = 17\)
        & \(p = 29\)
        & \(p = 41\)
        \\
        \hline
        \(\{K_{1}\}\) 
        & \(\subgrp{\tau,\zeta}\)
        & \TypeEzerosquared
        & \TypeEzerosquared
        & \TypeEzerosquared
        & \TypeEzerosquared
        & \TypeEzerosquared
        \\
        \(\{K_{2}, K_{4}, K_{6}\}\)
        & \(\subgrp{\sigma,\tau}\)
        & (loops)
        & (loops)
        & (loops)
        & (loops)
        & (loops)
        \\
        \(\{K_{3}, K_{5}, K_{7}\}\) 
        & \(\subgrp{\sigma,\tau}\)
        & \TypeEsquared
        & \TypeEtwelvecubedsquared
        & \TypeEsquared
        & \TypeEsquared
        & \TypeEsquared
        \\
        \(\{K_{8}, K_{9}\}\)
        & \(\subgrp{\tau\zeta, \zeta^2}\)
        & \TypeIV
        & \TypeIV
        & \TypeIV
        & \TypeVI
        & \TypeIV
        \\
        \(\{K_{i} : 10\le i\le 15\}\)
        & \(\subgrp{\sigma\tau}\)
        & \TypeI
        & \TypeIV
        & (loops)
        & \TypeI
        & \TypeIII
        \\
        \hline
    \end{tabular}
    \caption{Edge data for the \TypeV vertex}
    \label{tab:TypeV}
\end{table}

\begin{remark}
    To see the \TypeV neighbourhood in Figure~\ref{fig:TypeV}
    as a specialization of the \TypeIV diagram (Figure~\ref{fig:TypeIV}):
    \begin{itemize}
        \item
            the \TypeExE neighbour specializes to \(\classof{\EC}^2\),
            where \(\EC\) has \(j\)-invariant 54000
            and an endomorphism of degree 3;
        \item
            one of the \TypeIV neighbours degenerates to \TypeEzerosquared;
        \item
            the other two \TypeIV neighbours merge, yielding a weight-2 edge;
        \item
            one of the \TypeI neighbours specializes
            to \TypeV, yielding a loop;
        \item
            the other two \TypeI neighbours merge,
            yielding a weight-6 edge.
    \end{itemize}
\end{remark}

\subsection{The \texorpdfstring{\TypeVI}{Type-VI} vertex}

The unique \TypeVI vertex is \(\classof{\Jac{\XC_{VI}}}\),
where \(\XC_{VI} = \XC_{IV}(v_{VI})\)
with \(v_{VI} = (\zeta_{12}^2 + \zeta_{12} + 1)/\sqrt{2}\)
where \(\zeta_{12}^2 = \zeta_{6}\).
This curve is isomorphic to Bolza's \TypeVI normal form \(y^2 = x(x^4+1)\).
We have \(\RAut(\Jac{\XC_{VI}}) = \subgrp{\sigma,\rho,\omega} \cong S_4\),
where \(\sigma\) and \(\rho\) are inherited from \(\XC_{III}\)
and \(\omega\) is an order-4 automorphism acting as 
\[
    \omega_*: 
    x 
    \longmapsto
    (x - (\sqrt{2}+1))
    /
    ((\sqrt{2} - 1)x + 1)
\]
on \(x\)-coordinates.
Specializing the splittings of~\S\ref{sec:TypeI}
at \(s = s_{IV}(v_{VI}) = -\zeta_{12}^3\sqrt{2} - \zeta_{12}^3\)
and \(t = t_{IV}(v_{VI}) = 2\sqrt{2}+3\),
we see that \(\RAut(\Jac{\XC_{VI}})\) acts as
\begin{align*}
    \sigma_* & = (1)(2)(3)(4)(5)(6)(7)(8,9)(10,11)(12,13)(14,15)
    \,,
    \\
    \rho_*   & = (1,9,8)(2,15,14)(3)(4,12,13)(5)(6,10,11)(7)
    \,,
    \\
    \omega_* & = (1,4)(2,14,7,15)(3,10,6,11)(5)(8,9,13,12)
\end{align*}
on kernel indices. 
Table~\ref{tab:TypeVI} describes the kernel orbits and edges.
It is interesting to compare this with Table~\ref{tab:TypeIV},
to see how the various neighbours degenerate, specialize, and merge.
The \TypeEsquared neighbour is special:
it is \(\classof{\EC^2}\) where \(\EC\) is an elliptic curve of
\(j\)-invariant \(8000\);
it is \(\Phi\),
because \(\EC\) has a degree-3 endomorphism.
Pushing one step beyond the \TypeIV neighbours,
we find new \TypeI and \TypePhilabel{} vertices connected to
\(\classof{\EC^2}\),
and we thus complete the neighbourhood 
shown in Figure~\ref{fig:TypeVI}.

\begin{table}[ht]
    \centering
    \begin{tabular}{c|c|c@{\;}c@{\;}c}
        \multirow{2}{*}{Kernel orbit}
        & Stabilizer
        & \multicolumn{3}{c}{Codomain type}
        \\
        & (conjugate)
        & General \(p\) & \(p = 7\) & \(p = 13, 29\)
        \\
        \hline
        \(\{K_{1}, K_{4}, K_{8}, K_{9}, K_{12}, K_{13}\}\) 
        & \(\subgrp{\sigma,\omega^2}\)
        & \TypeEsquared
        & \TypeEtwelvecubedsquared
        & \TypeEsquared
        \\
        \(\{K_{3}, K_{6}, K_{10}, K_{11}\}\)
        & \(\subgrp{\sigma,\rho}\)
        & \TypeIV
        & (loops)
        & \TypeIV/\textsf{V}
        \\
        \(\{K_{2}, K_{7}, K_{14}, K_{15}\}\)
        & \(\subgrp{\sigma,\rho}\)
        & \TypeIV
        & (loops)
        & \TypeV/\textsf{IV}
        \\
        \(\{K_{5}\}\)
        & \(\RAut(\Jac{\XC_{VI}})\)
        & (loop)
        & (loop)
        & (loop)
    \end{tabular}
    \caption{%
        Edge data for the \TypeVI vertex.
        For \(p = 13\) and \(p = 29\),
        one of the two \TypeIV
        neighbours specializes to \TypeV,
        depending on the choice of \(\zeta_{12}\).
    }
    \label{tab:TypeVI}
\end{table}

\begin{figure}[ht]
    \centering
    \begin{tikzpicture}[
            >={angle 60},
            thick,
            vertex/.style = {circle, draw, fill=white, inner sep=0.5mm, minimum size=7mm},
            wt/.style = {fill=white, anchor=center, pos=0.5, minimum size=1mm, inner sep=1pt}
        ]
        \node (s) [
            regular polygon,
            regular polygon sides=6,
            minimum size=52mm,
            above
        ] at (0,0) {};
        \node (VI) [vertex] at (s.corner 3) {$VI$};
        \node (eephi) [vertex, dotted] at (s.corner 6) {\TypeEsquaredlabel \TypePhilabel};
        \node (IVa) [vertex, dotted] at (s.corner 2) {$IV$};
        \node (IVb) [vertex, dotted] at (s.corner 4) {$IV$};
        \node (I) [vertex, dotted] at (s.corner 5) {$I$};
        \node (newphi) [vertex, dotted] at (s.corner 1) {\TypePhilabel} ;

        \draw[->] (VI) edge[out=150,in=210,loop, looseness=6] (VI) ;
        \draw[->] (eephi) edge[out=330,in=30,loop, looseness=6] (eephi) ;

        \draw[->] (VI) edge[bend right=-24] node[wt]{4} (IVa) ;
        \draw[->] (IVa) edge[bend right=8, dotted] (VI) ;

        \draw[->] (VI) edge[bend right=24] node[wt]{4} (IVb) ;
        \draw[->] (IVb) edge[bend right=-8, dotted] (VI) ;

        \draw[->] (IVa) edge[bend right=-24, dotted] node[wt]{3} (newphi) ;
        \draw[->] (newphi) edge[bend right=8, dotted] (IVa) ;

        \draw[->] (IVa) edge[bend right=36, dotted] node[wt]{3} (I) ;
        \draw[->] (I) edge[bend right=-18, dotted] (IVa) ;

        \draw[->] (IVb) edge[bend right=24, dotted] node[wt]{3} (I) ;
        \draw[->] (I) edge[bend right=-8, dotted] (IVb) ;

        \draw[->] (IVb) edge[bend right=36, dotted] node[wt]{3} (newphi) ;
        \draw[->] (newphi) edge[bend right=-18, dotted] (IVb) ;

        \draw[->] (eephi) edge[bend right=8, dotted] node[wt]{2} (newphi) ;
        \draw[->] (newphi) edge[bend right=-24, dotted] (eephi) ;

        \draw[->] (eephi) edge[bend right=-8, dotted] node[wt]{2} (I) ;
        \draw[->] (I) edge[bend right=24, dotted] (eephi) ;

        \draw[->] (VI) edge[bend right=-8] node[wt]{6} (eephi) ;
        \draw[->] (eephi) edge[bend right=-8, dotted] (VI) ;

    \end{tikzpicture}
    \caption{The neighbourhood of the \TypeVI vertex.
        The dotted neighbours change type for \(p = 7\), \(13\), and \(29\)
        (see Table~\ref{tab:TypeVI}).
    }
    \label{fig:TypeVI}
\end{figure}

\subsection{The \texorpdfstring{\TypeII}{Type-II} vertex}

The unique \TypeII vertex
is \(\classof{\Jac{\XC_{II}}}\)
where \(\XC_{II}\) is defined  by \(y^2 = x^5 - 1\);
we have
\(\RAut(\Jac{\XC_{II}}) = \subgrp{\zeta} \cong C_5\),
where \(\zeta\) acts as
\[
    \zeta_*: x \longmapsto \zeta_5 x
    \,.
\]
The 15 kernel subgroups of \(\Jac{\XC_{II}}[2]\)
form three orbits of five under the action of \(\zeta\).
We fix orbit representatives
\begin{align*}
    K_1 & = 
    \{ x-1, (x-\zeta_5)(x-\zeta_5^2), (x - \zeta_5^3)(x - \zeta_5^4)\}
    \,,
    \\
    K_2 & = 
    \{ x-1, (x-\zeta_5)(x-\zeta_5^3), (x - \zeta_5^2)(x - \zeta_5^4)\}
    \,,
    \\
    K_3 & = 
    \{ x-1, (x-\zeta_5)(x-\zeta_5^4), (x - \zeta_5^2)(x - \zeta_5^3)\}
    \,,
\end{align*}
and let \(\phi_{i}: \Jac{\XC_{II}} \to \AV_{i} := \Jac{\XC_{II}}/K_i\)
be the quotient isogenies for \(1 \le i \le 3\).
Equation~\eqref{eq:ratio-principle} tells us that 
\(w(\classof{\phi_{i}}) = 5\) for each \(i\).

\begin{figure}[ht]
    \centering
    \begin{tikzpicture}[
            >={angle 60},
            thick,
            vertex/.style = {circle, draw, fill=white, inner sep=0.5mm, minimum size=7mm},
            wt/.style = {fill=white, anchor=center, pos=0.5, minimum size=1mm, inner sep=1pt}
        ]
        \node (s) [
            regular polygon,
            regular polygon sides=3,
            minimum size=42mm,
            above
        ] at (0,0) {};
        \node (c) [vertex] at (s.center) {$II$};
        \foreach \i in {1,2,3}{
            \node (s\i) [vertex, dotted] at (s.corner \i) {$A$};
            \draw[->] (c) edge[bend right=12] node[wt]{5} (s\i) ;
            \draw[->] (s\i) edge[bend right=12, dotted] (c) ;
        }
    \end{tikzpicture}
    \caption{%
        The neighbourhood of the (unique) \TypeII vertex.
    }
    \label{fig:TypeII}
\end{figure}
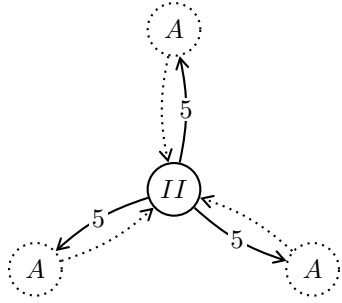

The neighbourhood of \(\classof{\Jac{\XC_{II}}}\)
is shown in Figure~\ref{fig:TypeII}.
Generally, the \(\classof{\AV_{i}}\) are \TypeA
(because the stabilizer of each orbit is trivial),
but for \(p = 19\), \(29\), \(59\), \(79\), and \(89\)
the codomain types change (see Table~\ref{tab:TypeII-codomains}).
Note that at \(p = 19\),
the codomain \(\AV_{2}\) becomes isomorphic to \(\AV_{0}\),
so \(\classof{\phi_{2}}\) becomes a weight-5 loop.

\begin{table}[ht]
    \centering
    \begin{tabular}{r|cccccc}
        {Characteristic \(p\)} & \(19\) & \(29\) & \(59\) & \(79\) & \(89\) & Other
        \\
        \hline
        \hline
        Type of \(\classof{\AV_{1}}\)
        &
        \TypeI   & \TypeI & \TypeI & \TypeI & \TypeA & \TypeA
        \\
        \hline
        Type of \(\classof{\AV_{2}}\)
        &
        \TypeII  & \TypeI & \TypeA & \TypeA & \TypeI & \TypeA
        \\
        \hline
        Type of \(\classof{\AV_{3}}\)
        &
        \TypeIII & \TypeA & \TypeI & \TypeA & \TypeA & \TypeA
        \\
        \hline
    \end{tabular}
    \caption{Types of neighbours of the \TypeII vertex.}
    \label{tab:TypeII-codomains}
\end{table}

\bibliographystyle{plain}
\bibliography{references}

\end{document}